\documentclass[review]{elsarticle}

\usepackage{lineno,hyperref}

\usepackage{color}

\usepackage{multirow}

\usepackage{mathtools}

\usepackage{amsmath}
\usepackage{amssymb}
\usepackage{graphicx}
\usepackage{algorithm}
\usepackage{psfrag}
\usepackage{subfig}
\usepackage{fullpage}

\journal{Journal of Computational Physics}

\begin{document}

\begin{frontmatter}

\title{Extending fields in a level set method by solving a biharmonic equation}

\author{Timothy J. Moroney}\corref{cor1}
\cortext[cor1]{Corresponding author}
\ead{t.moroney@qut.edu.au}
\author{Dylan R. Lusmore\corref{dummy}}
\author{Scott W. McCue\corref{dummy}}
\author{D.L. Sean McElwain\corref{dummy}}
\address{School of Mathematical Sciences, Queensland University of Technology (QUT), Brisbane QLD 4001, Australia}

\begin{abstract}
We present an approach for computing extensions of velocities or other fields in level set methods by solving a biharmonic equation.  The approach differs from other commonly used approaches to velocity extension because it deals with the interface fully implicitly through the level set function.  No explicit properties of the interface, such as its location or the velocity on the interface, are required in computing the extension.  These features lead to a particularly simple implementation using either a sparse direct solver or a matrix-free conjugate gradient solver.  Furthermore, we propose a fast Poisson preconditioner that can be used to accelerate the convergence of the latter.

We demonstrate the biharmonic extension on a number of test problems that serve to illustrate its effectiveness at producing smooth and accurate extensions near interfaces.  A further feature of the method is the natural way in which it deals with symmetry and periodicity, ensuring through its construction that the extension field also respects these symmetries.
\end{abstract}


\end{frontmatter}


\section{Introduction}
Moving boundary problems, or interface problems, arise in the mathematical modelling of many interesting natural phenomena, including Hele-Shaw flows \cite{paterson1981radial, hou1994removing, li2009control, dallaston2013bubble}, multiphase fluid flows \cite{hou2001boundary}, thin film flows \cite{bertozzi1997linear, diez2001contact, mayo2013gravity}, melting and solidification problems \cite{langer1980instabilities, chen1997simple, gibou2003level, back2014effect}, melanoma growth \cite{amar2011contour} and many others.

\begin{figure}
\centering
\begin{tabular}{cc}
\subfloat[Field known inside, requires extending out]{\includegraphics[trim = 0cm 2cm 0cm 2cm, clip, width=0.4\linewidth]{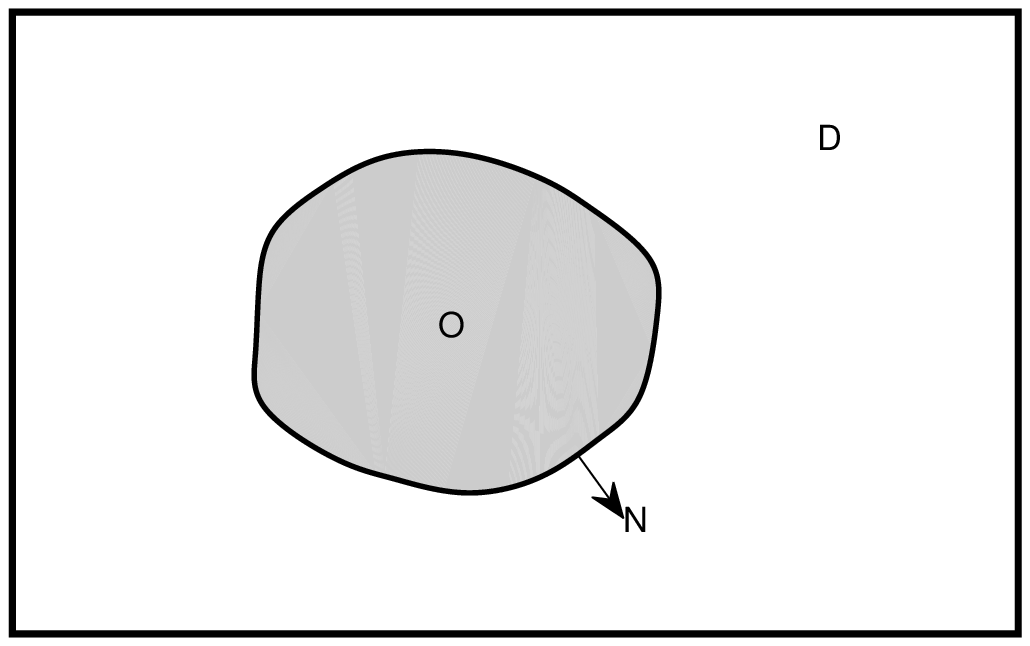}} &
\subfloat[Field known outside, requires extending in]{\includegraphics[trim = 0cm 2cm 0cm 2cm, clip, width=0.4\linewidth]{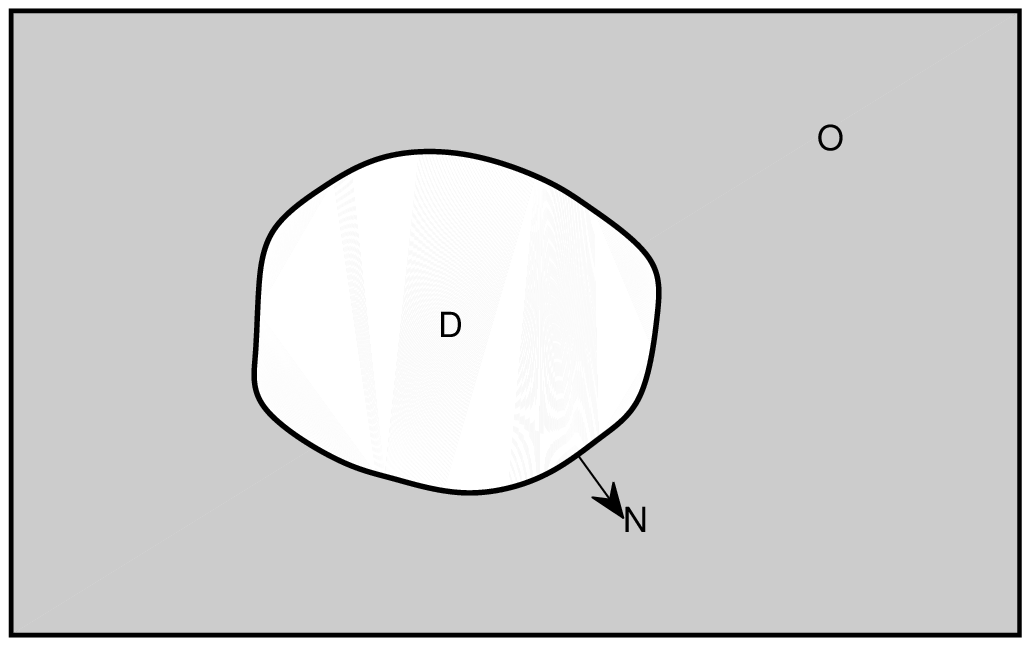}}
\end{tabular}
\caption{Schematics of typical moving boundary problems.  The interface moves with normal velocity $V_n = \mathbf{F} \cdot \mathbf{n}$.  The velocity field $\mathbf{F}$ is calculated in the shaded region $\Omega$ from the solution to a field equation.  The level set method then requires $\mathbf{F}$ to be smoothly extended to all of $D$.  Here two common scenarios are illustrated, requiring the field to be extended either ``out'' (a) or ``in'' (b).}
\label{fig:schematic}
\end{figure}

Schematics of two particular configurations of a moving boundary problem on a two-dimensional domain $D$ are shown in Figure \ref{fig:schematic}.  We suppose that $D$ is partitioned into two regions, and that the interface $\partial \Omega$ between the two evolves under a normal velocity given by
\begin{equation}
\label{eq:normal_velocity}
V_n = \mathbf{F}\cdot\mathbf{n} \quad \textrm{on } \partial \Omega\,.
\end{equation}
In the shaded region $\Omega$ the velocity field $\mathbf{F}$ is calculated from the solution of a field equation, which we suppose is a parabolic or elliptic partial differential equation.  (We do not consider hyperbolic equations in this paper, since as noted by Gibou et al. \cite{Gibou2013} they require further considerations to ensure that conservation properties are conserved, in order that the correct speed of propagation is calculated in the presence of shocks.)  Having calculated the velocity field $\mathbf{F}$, the normal velocity of the interface $V_n$ is then determined from \eqref{eq:normal_velocity}.  In this way, the interface evolves over time and, depending on the phenomena under consideration, may develop fingering instabilities, splitting, merging, and other interesting features.  Accurately tracking the evolution of the interface over time is a key requirement of any numerical solution method for such a problem.

\psfrag{N}{$\mathbf{n}$}
\psfrag{O}{$\Omega$}
\psfrag{D}{$\!\!\!\!D \backslash \Omega$}

The level set method is one numerical method that has proved very popular for solving interface problems \cite{osher1988fronts, Osher2006LevelSet,Sethian1999LevelSet}.  Rather than developing an explicit representation for the interface location, the level set method represents the interface \emph{implicitly}, as the zero level set of an auxiliary function $\phi$ -- the so-called level set function.  This approach has the distinct advantage that the location of the interface, as well as splitting, merging, and other such complex behaviours as the interface evolves, are all handled naturally by the level set function discretised on a standard rectangular grid.  The evolution of $\phi$ under the external velocity field $\mathbf{F}$ is governed by the level set equation
\begin{equation}
\label{eq:level_set_equation}
\frac{\partial \phi}{\partial t} + \mathbf{F} \cdot \nabla \phi = 0\,.
\end{equation}
A typical numerical scheme proceeds by solving the field equation to determine the velocity field, followed by updating the level set according to \eqref{eq:level_set_equation}, discretised using finite differences.  There is, however, a catch.  Although the interfacial velocity $V_n$ has physical significance only on the interface itself, the level set equation \eqref{eq:level_set_equation} requires a field $\mathbf{F}$ defined on all $D$.  Hence, a means of extrapolating, or extending, the velocity field $\mathbf{F}$ from $\Omega$ to $D$ is required to complete the level set implementation.  This extension velocity need not have any physical significance itself, and is required only to be a smooth, well-defined field on all of $D$ that accurately represents the true velocity in a small neighbourhood of the interface.

Many numerical schemes have been proposed for constructing extension velocities.
One approach is to solve the equation
\[
\nabla f \cdot \nabla \phi = 0
\]
by a process known as the Fast Marching Method \cite{Adalsteinsson1999FastConstruction}, which sweeps forward in the upwind direction from the interface, building the extension within an evolving narrow band.  This approach is motivated by the desire to maintain the level set $\phi$ as a signed distance function, satisfying $\| \nabla \phi\| = 1$, under evolution by \eqref{eq:level_set_equation}.  Higher order variants of the method can also be constructed provided appropriate care is taken with the direction of propagation near the interface for the higher order normal derivatives \cite{McCaslin2014FastMarching}.

An alternative approach to constructing extension velocities is to solve an advection equation
\begin{equation}
\label{eq:advection_extension}
\frac{\partial f}{\partial \tau} + H(\phi) \nabla f \cdot \mathbf{n} = 0
\end{equation}
to steady-state in pseudo-time $\tau$, 
where $H$ denotes the Heaviside function \cite{fedkiw1999non}.  
The resulting field is constant along the normal directions, carrying the interfacial value along these characteristics.  Higher order schemes can be built upon this foundation by a sequence of extrapolations, first with higher order normal derivatives, through lower order normal derivatives and finally the value itself \cite{Aslam2004PartialDifferentialEquation}.  Alternatively, the steady version of \eqref{eq:advection_extension} may be solved directly
using a fast sweeping method that carefully chooses and adjusts the required upwind stencils and systematically alternates through different node orderings to cover all characteristic directions \cite{Aslam2014StaticPDE}.

A common feature of these approaches is the need to determine certain explicit properties related to the interface, such as its location, the location of the nearest node to it, or the value of the field on the interface.

This is in contrast with the level set method itself, which deals with the interface in an entirely implicit fashion. 
Hence we are motivated to pursue a velocity extension approach that also deals with the interface implicitly, and does not rely on an explicit representation of the interface.

A further point we address concerns problems possessing symmetry or periodicity.  For these problems, it is natural to exploit this symmetry in the implementation of the level set method by solving over a reduced computational domain.  This approach is valid provided the extension field $\mathbf{F}$ \emph{also} possesses the required symmetry.  Hence, it is an important consideration when constructing velocity extensions for such problems that they conform to the prescribed symmetry or periodicity conditions on the domain boundary.

In this paper we present an extension method that deals only with the interface implicitly, making for a scheme that is straightforward to implement and fully in keeping with the spirit of the level set method.  Furthermore, the extension field naturally conforms to any desired symmetry or periodicity conditions on the boundary of the computational domain.  The method should be useful for extending velocity fields, but also temperatures, mass fractions, or other fields relevant to the application.  In the next section we outline the method, introducing the theoretical background, and we also propose a simple numerical implementation.  In Section \ref{sec:numerical} we present five examples that illustrate the effectiveness of the method.  Finally, we draw our conclusions in Section \ref{sec:conclusion}.

\section{Biharmonic extension}
\subsection{Formulation of our approach}
\label{sec:biharmonic}
In the problem as described thus far, values of some vector field $\mathbf{F}$ are known on one side of an interface and must be extended in a sensible way to the entire domain.  From this point onwards, we consider only a scalar field, which may be a component of $\mathbf{F}$, or more generally any scalar field for which an extension is required.  Values of this field are known at all nodes of the finite difference grid in $\Omega$.  The extension problem is then: Given known values $f_k$ for $k=1,\ldots, N$ at nodes $\mathbf{x}_k \in \Omega$, construct a suitably smooth function $f: D \to \mathbb{R}$ such that
\begin{equation}
\label{eq:interpolation_constraints}
f(\mathbf{x}_k) = f_k,\qquad \textrm{for } k=1,\ldots, N\,.
\end{equation}

To motivate our approach, we begin by considering the simpler problem of \emph{interpolation}, such as would be the case if the domain was partitioned as in Figure \ref{fig:schematic}(b), that is, with the field values $f_k$ known on the outside of $\partial \Omega$.  The problem of determining field values on the inside is then an interpolation, since the interpolation points lie within the convex hull of the prescribed data.  In two dimensions, a natural choice of interpolant is the thin plate spline \cite{Bookstein1989PrincipalWarps}.  Formulated in terms of radial basis functions (RBFs) \cite{Buhmann2000RadialBasis}, the thin plate spline interpolant takes the form
\begin{equation}
\label{eq:radial_basis_function}
f(\mathbf{x}) = \sum_{k=1}^N c_k \, \Phi(\|\mathbf{x} - \mathbf{x}_k\|)\,,
\end{equation}
where $\Phi: \mathbb{R}^+ \to \mathbb{R}$, a function of the distance between points, is given by
\begin{equation}
\label{eq:thin_plate_spline}
\Phi(r) = r^2 \log(r)
\end{equation}
and the $N$ coefficients $c_k$ are determined by enforcing the $N$ interpolation constraints \eqref{eq:interpolation_constraints}.

The terminology ``thin plate spline'' derives from one application of this interpolant, which when augmented by an affine polynomial term (i.e. $f(\mathbf{x}) + a_0 + b_1 x + b_2 y$)
describes the shape of a thin metal plate, infinite in extent, which is constrained to pass through the set of interpolation points.  In particular, the function minimises the ``bending energy''
\begin{equation}
\iint_{\mathbb{R}^2} (f_{xx}^2 + 2f_{xy}^2 + f_{yy}^2)\, \textrm{d}x\, \textrm{d}y
\end{equation}
over all possible interpolants $f$ \cite{Bookstein1989PrincipalWarps}.  

The thin plate spline interpolant can in principle be determined directly by solving an $N \times N$ linear system for the coefficients $c_k$ in \eqref{eq:radial_basis_function} (or, if the affine polynomial term is included, an $(N+3) \times (N+3)$ linear system, which allows for the imposition of three further constraints that force the solution to be asymptotically flat as $\|\mathbf{x}\| \to \infty$).  There are, however, several limitations that render this approach impractical for the present application of computing extension velocities.  First, it is well known that the conditioning of the RBF linear system is typically very poor, and increasingly so for larger $N$.  Second, whenever the RBF has global support (as it does for the thin plate spline), the coefficient matrix is dense, and hence storage and factorisation costs are prohibitive for large $N$.  Third, in general, computing an extension velocity requires \emph{extrapolation} rather than interpolation, such as when the domain is partitioned as in Figure \ref{fig:schematic}(a), and the thin plate spline is unbounded away from the prescribed data.

These limitations may be overcome by considering the thin plate spline interpolant from an alternative point of view \cite{Gaspar2000MultilevelBiharmonic}.  The thin plate spline RBF \eqref{eq:thin_plate_spline} is (to within a multiplicative constant) a fundamental solution of the biharmonic equation $\nabla^4 \Phi = 0$ in $\mathbb{R}^2$.  Hence, by superposition, the linear combination \eqref{eq:radial_basis_function} also solves the biharmonic equation in $\mathbb{R}^2$, except possibly at the points $\mathbf{x}_k$ themselves.  Being a function defined on all of $\mathbb{R}^2$, it satisfies no particular boundary conditions and indeed as already noted, it is not bounded for large $\|\mathbf{x}\|$.  However, if we instead consider the problem on the finite domain $D \subset \mathbb{R}^2$, we may now enforce boundary conditions on $f$, which can be used to ``clamp'' the thin plate spline, for example by enforcing that it vanishes on $\partial D$.

The idea then, is to determine the extension field $f$ as the solution to the biharmonic equation on $D \backslash \{\mathbf{x}_1, \ldots, \mathbf{x}_N\}$, subject to suitable boundary conditions on $\partial D$ and also subject to the interpolation constraints \eqref{eq:interpolation_constraints}.  We shall call this function $f$ the \emph{biharmonic extension} which, supposing for the present time, that homogeneous Dirichlet boundary conditions are imposed, satisfies
\begin{equation}
\label{eq:biharmonic}
\nabla^4 f = 0 \qquad \textrm{in } D \backslash \{\mathbf{x}_1, \ldots, \mathbf{x}_N\}
\end{equation}
subject to
\begin{equation}
\label{eq:biharmonic_boundary_conditions}
f = 0 \quad \textrm{and} \quad \nabla^2 f = 0 \qquad \textrm{on } \partial D
\end{equation}
and
\begin{equation}
\label{eq:biharmonic_interpolation_constraints}
f(\mathbf{x}_k) = f_k,\qquad \textrm{for } k=1,\ldots, N\,.
\end{equation}
We note that this fourth order problem is well-posed with the specification of both \eqref{eq:biharmonic_boundary_conditions} and \eqref{eq:biharmonic_interpolation_constraints}.  Indeed, the problem could instead be posed by surrounding each interpolation point $\mathbf{x}_k$ with a circle $\Gamma_k$ of radius $\epsilon$, and imposing
\begin{equation}
\label{eq:biharmonic_interpolation_boundary_conditions}
f = f_k \quad \textrm{and} \quad  \partial f / \partial n = 0 \qquad \textrm{on } \Gamma_k
\end{equation}
instead of \eqref{eq:biharmonic_interpolation_constraints}.  The unique biharmonic extension would then be obtained in the limit $\epsilon \to 0$ \cite{Gaspar2000MultilevelBiharmonic}.

Solving \eqref{eq:biharmonic} -- \eqref{eq:biharmonic_interpolation_constraints} to determine the biharmonic extension is very convenient in the context of the level set method.  The solution can be computed on the same finite difference grid already in place for the level set equation \eqref{eq:level_set_equation}.  Just as with the level set equation, there is no need to take special care with finite difference stencils near the interface.  Stencils can, and do, cross over the interface, and indeed this is the whole point of computing the extension -- it carries information from one side of the interface to the other.

Numerically, we simply solve the discrete biharmonic equation
\begin{align}
\label{eq:discrete_biharmonic}
\nonumber
&20f_{i,j} - 8f_{i-1,j} - 8f_{i,j+1} - 8f_{i,j-1} - 8f_{i+1,j} + 2f_{i+1,j+1} + 2f_{i-1,j-1} \\
&+ 2f_{i+1,j-1} + 2f_{i-1,j+1} +f_{i+2,j} + f_{i-2,j} + f_{i,j+2}  + f_{i,j-2}  = 0
\end{align}
at every node $(x_i, y_j)$ in $D \backslash \Omega$.  
Stencils near the domain boundary are modified according to the boundary conditions \eqref{eq:biharmonic_boundary_conditions}.  Where the stencil crosses the interface, the referenced values on the other side are known: they are the values $f_k$ that we wish to extend.  In this way, the biharmonic extension naturally incorporates two levels of information on one side of the interface (the values $f_{i-1, j}$ and $f_{i-2,j}$, say) and constructs a smooth extension that follows the trend of the data across the interface to the other side.  The resulting extension is, by construction, continuous in both value and derivative.

We also emphasise that at no time is the location of the interface, nor the value of the field on the interface, required in constructing the extension.  The biharmonic extension deals only with the field values on the level set grid, eliminating any need to explicitly locate, or map values to, the interface before extending across to the other side.
Indeed, a significant advantage of the biharmonic extension is its simplicity of implementation, as we now discuss.

\subsection{Efficient computation of biharmonic extension}
\label{subsec:iterative}
As already mentioned, forming the extension field by solving the biharmonic equation rather than with the RBF equation \eqref{eq:radial_basis_function} avoids the dense linear systems associated with RBF-based interpolation.  Instead, the coefficient matrix is the discrete biharmonic matrix arising from imposing \eqref{eq:discrete_biharmonic} at each node.  Nevertheless, it is important to consider how this sparse, symmetric positive definite linear system can be solved in order to build the extension field efficiently.

For problems where the number of nodes at which the extension is required is small compared to the total problem size (i.e. where $\Omega$ comprises most of the computational domain $D$), it is efficient to simply form and factorise the coefficient matrix using standard sparse algorithms (e.g. MATLAB's sparse \verb"chol" command).  For extensions involving a large number of nodes, the matrix fill-in incurred by factorisation may render this approach unviable, and alternative solution methods may be required.  Fast finite difference methods for biharmonic equations over irregular domains have received some attention in the literature, for example \cite{Mayo1984FastPoission, mayo1992fast, chen2008fast}.  The idea is to embed the irregular domain in a rectangular region and apply fast Poisson solvers over the rectangle (exploiting the fact that the biharmonic operator is the square of the Laplacian), with methods differing in the way the embedding is carried out.  Other approaches for building thin plate spline interpolants through the biharmonic equation have used adaptive quadtree and octtree grids \cite{Gaspar2000MultilevelBiharmonic}.  Further ideas can be found from the field of image ``inpainting'', where Poisson or biharmonic equations are often solved to ``fill in'' missing details from images.  Here, methods include multigrid methods \cite{Mainberger20111859} and methods based on conjugate gradients and pseudo-timestepping \cite{schmaltz2014understanding}.  Hence, there is great scope for further research into developing fast methods for computing these extension fields.

\begin{algorithm}
\caption{Complete MATLAB code to compute the biharmonic extension}
\label{al:biharmonic}
\begin{footnotesize}
\begin{verbatim}
% Biharmonically extend f to D from known values in Omega
% f:    m x n array with values in Omega filled, values in D \ Omega currently zero
% idx:  index of entries in f corresponding to nodes in D \ Omega
% tol:  iterative solve tolerance

b = -biharmonic_stencil(f); b = b(idx);                         % Set right hand side
f(idx) = pcg(@matvec, b, tol, max(m,n), @fast_poission);        % Solve with conjugate gradient

function B = biharmonic_stencil(f)
    i = 3:m+2; j = 3:n+2;
    F = zeros(m+4, n+4); F(i,j) = f;
    F([1 m+4],:) = -F([3 m+2],:);                               % Pad around boundaries
    F(:,[1 n+4]) = -F(:,[3 n+2]);

    i = 2:m+3; j = 2:n+3;
    L = 4*F(i,j) - F(i-1,j) - F(i+1,j) - F(i,j-1) - F(i,j+1);   % Laplacian stencil
    i = 2:m+1; j = 2:n+1;
    B = 4*L(i,j) - L(i-1,j) - L(i+1,j) - L(i,j-1) - L(i,j+1);   % And again
end

function x = fast_poission(b)
    Lm = 2*(1-cos((1:m)*pi/(m+1)));                             % Eigenvalues
    Ln = 2*(1-cos((1:n)*pi/(n+1)));
    LL = Lm'*ones(1,n) + ones(m,1)*Ln;

    B = zeros(m,n); B(idx) = b;                                 % Map to D

    Z = dst(dst(B')');                                          % Inversion
    Z = Z ./ LL.^2;
    Z = idst(idst(Z')');

    x = Z(idx);                                                 % Map back to D\Omega
end

function b = matvec(u)
    U = zeros(m,n); U(idx) = u;
    B = biharmonic_stencil(U);
    b = B(idx);
end
\end{verbatim}
\end{footnotesize}
\end{algorithm}

In the present context we have an appropriate rectangular region $D$, and so we opt for a method based on a simple embedding coupled with a fast Poisson solver.  The mapping from $D \backslash \Omega$ to $D$ is achieved by augmenting the right hand side vector of the linear system with zeros, while the mapping back to $D \backslash \Omega$ simply omits the additional entries in the vector.  The correspondence between problems on the two domains is hence only approximate, so our approach is one of \emph{preconditioning}: using a fast Poisson solver on the rectangular region $D$ as a preconditioner for a conjugate gradient solver on $D \backslash \Omega$.

While undoubtedly more sophisticated schemes are possible, this approach to computing the biharmonic extension has a simple implementation, with the entire solution process described by the MATLAB code in Algorithm \ref{al:biharmonic}. We note that there is no need to form the 13-banded discrete biharmonic matrix at all with this iterative approach: the left hand side of \eqref{eq:discrete_biharmonic} can simply be evaluated on demand.


The cost of applying the fast Poisson preconditioner is $O(N\log N)$ operations for a grid of $N = n \times n$ nodes in total.  In the next section, we report for each of our numerical experiments the number of conjugate gradient iterations required for convergence using this preconditioner.  Typically we found that the required number of iterations was a slowly growing function of the linear dimension $n$.

\section{Numerical experiments}
In this section we demonstrate and discuss the effectiveness of the biharmonic extension approach on five distinct examples, in two and three dimensions, including problems exhibiting symmetry and periodicity.

\label{sec:numerical}
\subsection{Example 1, two-dimensional peanut}
\begin{figure}
\centering
\begin{tabular}{cc}
\subfloat[Reference field]{\includegraphics[trim= 2.5cm 0cm 2.5cm 0cm, clip, width=0.45\linewidth]{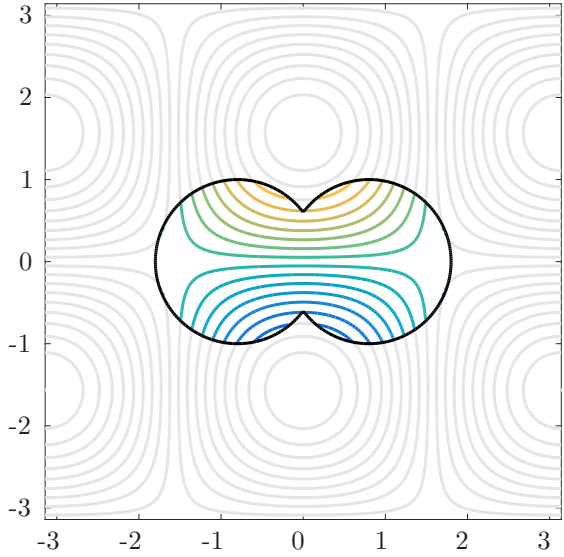}} &
\subfloat[Biharmonic extension subject to Dirichlet boundary conditions]{\includegraphics[trim= 2.5cm 0cm 2.5cm 0cm, clip, width=0.45\linewidth]{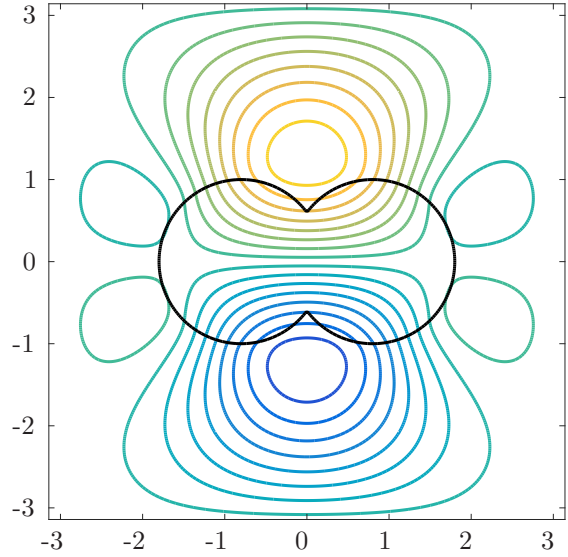}} \\
\subfloat[Biharmonic extension subject to Neumann conditions]{\includegraphics[trim= 2.5cm 0cm 2.5cm 0cm, clip, width=0.45\linewidth]{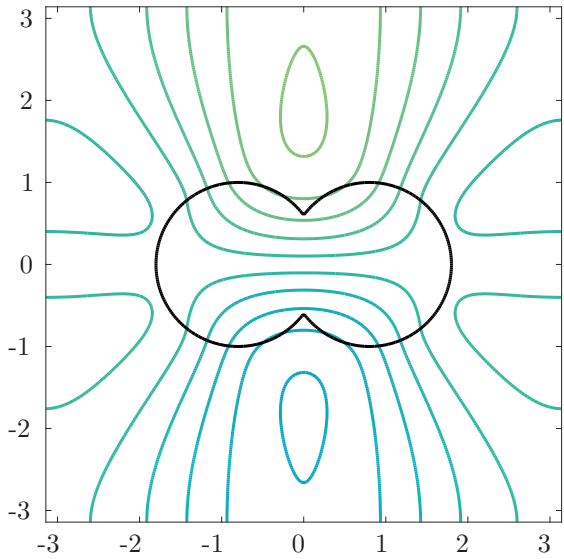}} &
\subfloat[Biharmonic extension subject to mixed Neumann and Dirichlet conditions]{\includegraphics[trim= 2.5cm 0cm 2.5cm 0cm, clip, width=0.45\linewidth]{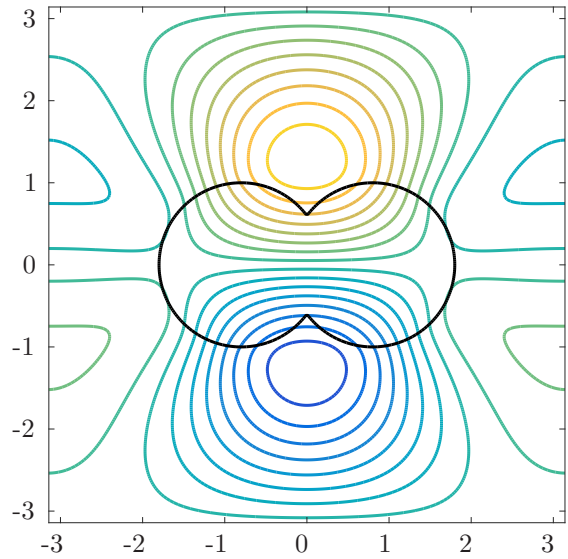}}\\
\end{tabular}
\caption{The reference field and biharmonic extension fields for Example 1 on a $256 \times 256$ grid, subject to Dirichlet, Neumann and mixed (Neumann on sides, Dirichlet top and bottom) boundary conditions.  Contour levels are from $[-1, 1]$.}
\label{fig:example1}
\end{figure}

For our first example, taken from \cite{Aslam2014StaticPDE}, we consider the domain $D = [-\pi, \pi] \times [-\pi, \pi]$ with field values
\begin{equation}
\label{eq:f_example1}
f(x,y) = \cos(x)\sin(y)
\end{equation}
prescribed in the region
\[
\Omega = \left \{(x,y) \in D \left| \min\left(\sqrt{(x-0.8)^2 + y^2} - 1, \sqrt{(x+0.8)^2 + y^2} - 1\right) < 0 \right. \right \}\,.
\]
This gives rise to the interface shown in black in Figure \ref{fig:example1}(a), with the field known on the inside.  The goal is to smoothly extend the field to all of $D$.  We also show in grey the extension given by formula \eqref{eq:f_example1} itself for reference.  We will benchmark our computed extensions by comparing against this reference field in a small neighbourhood of the interface.

\psfrag{-3}[cb]{-3}
\psfrag{-2}[cb]{-2}
\psfrag{-1}[cb]{-1}
\psfrag{0}[cb]{0}
\psfrag{1}[cb]{1}
\psfrag{2}[cb]{2}
\psfrag{3}[cb]{3}
\psfrag{0.5}[cb]{0.5}
\psfrag{-0.5}[cb]{-0.5}
\psfrag{1.5}[cb]{1.5}
\psfrag{-1.5}[bc]{-1.5}

\psfrag{n3}[cb]{ -3}
\psfrag{n4}[cb]{ -4}
\psfrag{n5}[cb]{ -5}
\psfrag{n6}[cb]{ -6}
\psfrag{n7}[cb]{ -7}
\psfrag{n8}[cb]{ -8}

We test three different choices of boundary conditions on $D$, namely Dirichlet, Neumann and a mixture of Neumann on the sides and Dirichlet on the top and bottom.  (Modifications to the code in Algorithm \ref{al:biharmonic} for Neumann conditions are straightforward, requiring just a change to the array padding in the biharmonic stencil, and changing the eigenvalues and fast transforms from sine transforms to cosine transforms.)
Figure \ref{fig:example1}(b)--(d) illustrates the biharmonic extension fields under these three choices of boundary conditions.  The extension field shown in Figure \ref{fig:example1}(d) comes closest to representing the reference field in Figure \ref{fig:example1}(a), owing to the chosen boundary conditions matching the actual behaviour of formula \eqref{eq:f_example1}.  But very often in practice no such boundary information will be available, and the computational boundary is simply placed sufficiently far from the interface, as is the usual practice with level set methods. Hence the essential focus is on the behaviour of the extension in the neighbourhood of the interface and, to visual accuracy, the fields in Figure \ref{fig:example1}(b)--(d) do equally well in this regard.

\begin{table}
\centering
\caption{Iteration counts, errors and estimated orders of convergence for the biharmonic extension over a sequence of grids for Example 1.  Dirichlet, Neumann and mixed boundary conditions are tested.  All errors are calculated with respect to the reference field \eqref{eq:f_example1} in a neighbourhood of 4 grid cells of the interface.}
\begin{tabular}{l|rrrr}
\hline
Grid& $128 \times 128$ & $256 \times 256$ & $512 \times 512$ & $1024 \times 1024$ \\
\hline
& \multicolumn{4}{c}{Biharmonic extension (Dirichlet)}\\
\hline
Iterations &         50 &        112 &        238 &        543 \\
Error &   6.28e-02 & 1.77e-02 & 4.71e-03 & 1.18e-03 \\
Est. order &        -- &       1.83 & 1.91 & 2.00  \\
\hline
& \multicolumn{4}{c}{Biharmonic extension (Neumann)}\\
\hline
Iterations &         50 &        117 &        253 &        558 \\
Error &   5.70e-02 & 1.57e-02 & 4.10e-03 & 1.02e-03  \\
Est. order &        -- &       1.86 & 1.94 & 2.00  \\
\hline
& \multicolumn{4}{c}{Biharmonic extension (mixed)}\\
\hline
Iterations &         53 &        120 &        257 &        579 \\
Error &   5.57e-02 & 1.55e-02 & 4.09e-03 & 1.02e-03 \\
Est. order &        -- &       1.84 & 1.92 & 2.00 \\
\hline
\end{tabular}
\label{tab:example1}
\end{table}

Table \ref{tab:example1} confirms this observation numerically, by tabulating the maximum difference between each extension and the reference field in a neighbourhood of four grid cells of the interface.  As the number of grid divisions is doubled in each direction, the errors for all three choices of boundary condition are seen to reduce by approximately a factor of four, which is consistent with second order spatial accuracy.

Using the iterative solution strategy given in Algorithm \ref{al:biharmonic}, the number of conjugate gradient iterations required to compute the field on an $n \times n$ grid is seen to grow slightly faster than linearly with $n$, requiring for example 257 and 579 iterations for the mixed case on $512 \times 512$ and $1024 \times 1024$ grids respectively.  All solutions were solved to iterative tolerance $10^{-6}$ in relative norm, which was sufficient for all tabulated errors to be unchanged with further iteration.

In fact, for this example the direct approach of simply forming and factorising the discrete biharmonic matrix proves to be more efficient than the iterative solver.  In Table \ref{tab:example1_errors} we report the runtimes for both direct and iterative methods for the biharmonic extension, observing that the direct method is much faster for this problem.  For example, with $n = 1024$ the direct solver requires 12.4 seconds compared to 131.0 seconds for the iterative solver.  For comparison purposes, we also report in Table \ref{tab:example1_errors} the errors and runtimes for computing constant, linear and quadratic extrapolations by solving \eqref{eq:advection_extension}, and its higher order variants, as described by Aslam \cite{Aslam2004PartialDifferentialEquation}.  As expected, these methods are observed to converge at first, second, and third order respectively to the reference field in a neighbourhood of the interface.  Concerning runtime, the biharmonic extension computed iteratively is as fast, or faster, than all the extrapolation methods, while the biharmonic extension computed directly is a clear order of magnitude faster to compute than the other methods.  We also report in Table \ref{tab:example1_errors} the condition numbers of the discrete biharmonic matrix, estimated using MATLAB's \verb"condest" function.  These grow at the expected rate of 16 as the number of divisions is doubled, and are very similar to the condition numbers obtained from the squares of the eigenvalues of the discrete Laplacian on a rectangle, as in Algorithm \ref{al:biharmonic}.

\begin{table}
\centering
\caption{Runtimes, condition numbers, errors and estimated orders of convergence over a sequence of grids for Example 1 for the biharmonic extension (Neumann conditions) and the constant, linear and quadratic extrapolations.  All errors are calculated with respect to the reference field \eqref{eq:f_example1} in a neighbourhood of 4 grid cells of the interface.}
\begin{tabular}{l|*{4}{r}}
\hline
Grid & $128\times128$ & $256\times256$ & $512\times512$ & $1024\times1024$ \\
\hline
& \multicolumn{4}{c}{Biharmonic extension (Neumann)} \\
\hline
Runtime (direct) & 0.1 & 0.4 & 2.3 & 12.4   \\
Runtime (iterative) & 0.1 & 0.9 & 13.2 & 131.0 \\
Cond. Est. & 2.4e+07 & 3.6e+08 & 5.7e+09 & 9.1e+11 \\
Error & 5.70e-02 & 1.57e-02 & 4.10e-03 & 1.02e-03  \\
Est. Order & -- & 1.86 & 1.94 & 2.00  \\
\hline
& \multicolumn{4}{c}{Constant extrapolation} \\
\hline
Runtime & 0.2 & 1.0 & 13.0 & 128.5 \\
Error & 1.10e-01 & 6.06e-02 & 3.14e-02 & 1.57e-02  \\
Est. Order & -- & 0.86 & 0.95 & 1.00 \\
\hline
& \multicolumn{4}{c}{Linear extrapolation} \\
\hline
Runtime & 0.3 & 2.0 & 25.8 & 258.9 \\
Error & 2.09e-02 & 5.21e-03 & 1.34e-03 & 3.38e-04  \\
Est. Order & -- & 2.00 & 1.96 & 1.99   \\
\hline
& \multicolumn{4}{c}{Quadratic extrapolation} \\
\hline
Runtime & 0.5 & 3.0 & 40.8 & 392.5  \\
Error & 3.19e-03 & 4.58e-04 & 5.41e-05 & 6.73e-06  \\
Est. Order & -- & 2.80 & 3.08 & 3.01   \\
\hline
\end{tabular}
\label{tab:example1_errors}
\end{table}

\begin{figure}
\centering
\begin{tabular}{cc}
\subfloat[Biharmonic extension (Neumann)]{\includegraphics[trim= 2.5cm 0cm 2.5cm 0cm, clip, width=0.45\linewidth]{Example1_neumann_2.eps}} &
\subfloat[Constant extrapolation]{\includegraphics[trim= 2.5cm 0cm 2.5cm 0cm, clip, width=0.45\linewidth]{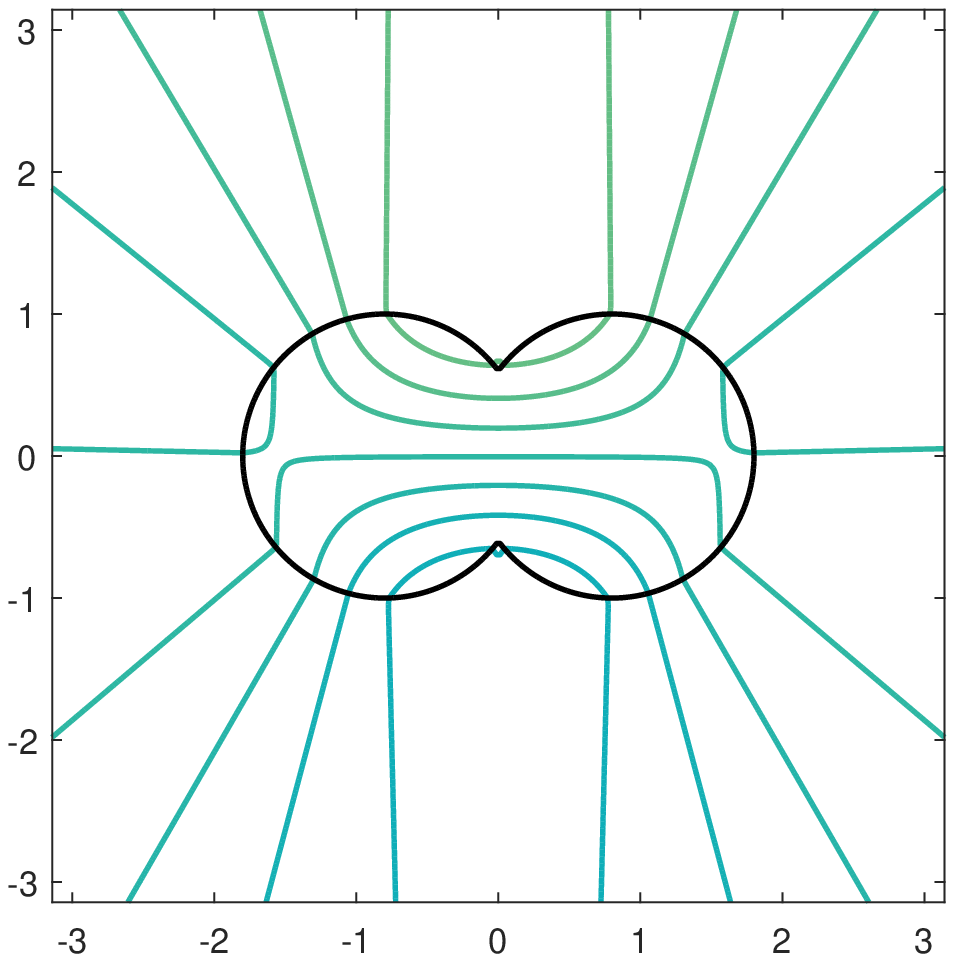}} \\
\subfloat[Linear extrapolation]{\includegraphics[trim= 2.5cm 0cm 2.5cm 0cm, clip, width=0.45\linewidth]{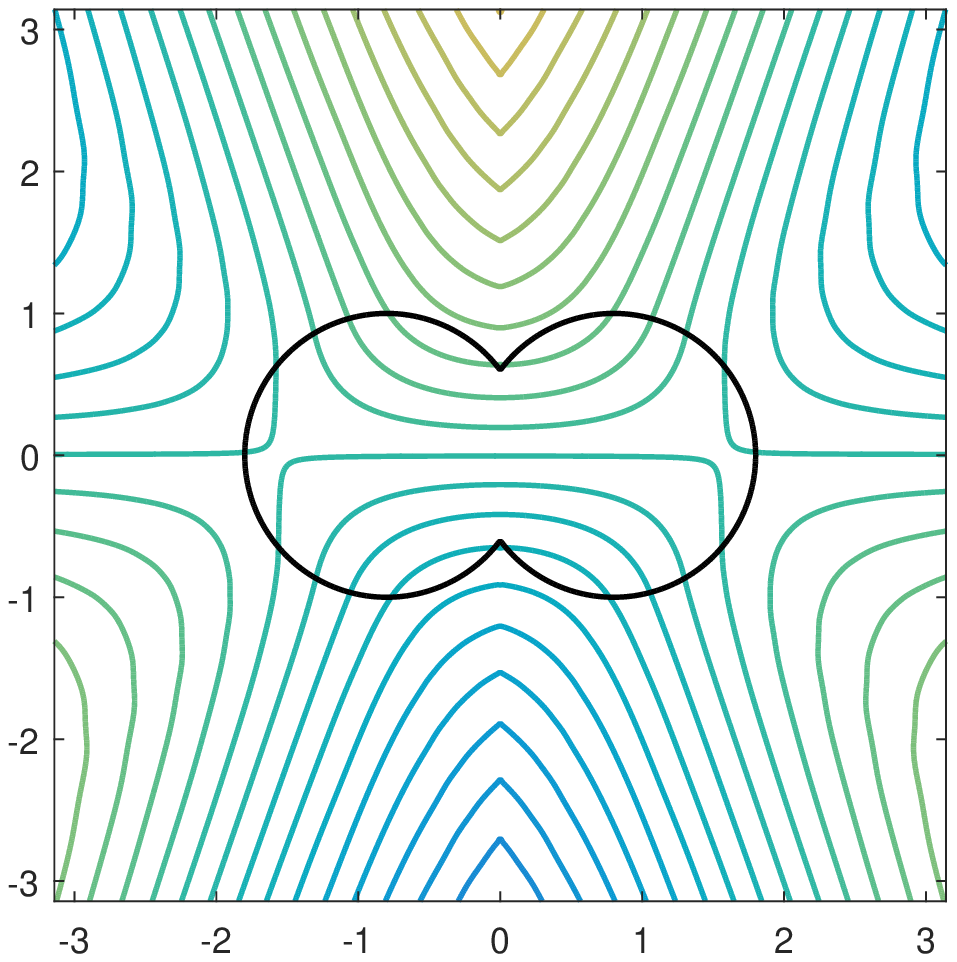}} &
\subfloat[Quadratic extrapolation]{\includegraphics[trim= 2.5cm 0cm 2.5cm 0cm, clip, width=0.45\linewidth]{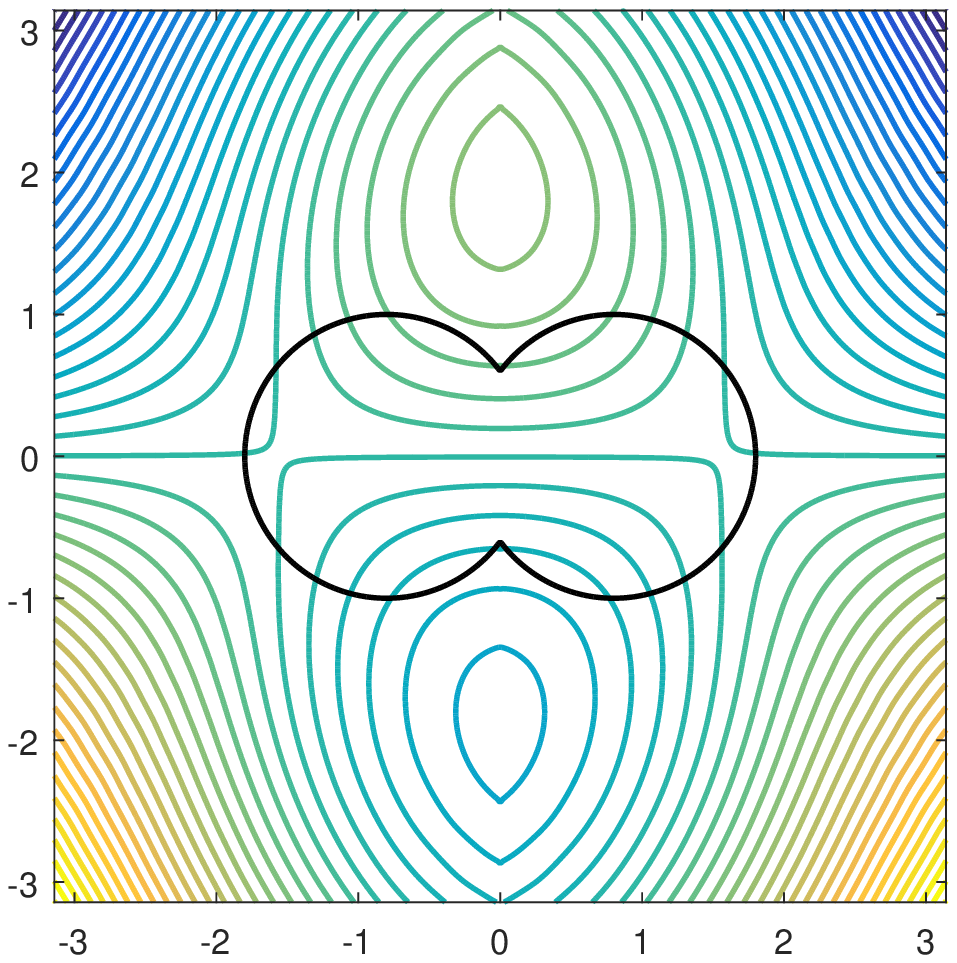}}\\
\end{tabular}
\caption{The biharmonic extension field for Example 1 on a $256 \times 256$ grid subject to Neumann boundary conditions, compared with constant extrapolation, linear extrapolation and quadratic extrapolation using the method of Aslam \cite{Aslam2004PartialDifferentialEquation}. Contour levels are now from $[-4, 4]$ for ease of comparison with \cite{Aslam2014StaticPDE} Fig. 3.}
\label{fig:example1b}
\end{figure}

The constant, linear and quadratic extrapolations are compared with the biharmonic extension (with Neumann conditions) in Figure \ref{fig:example1b}, which reproduces the results presented in \cite{Aslam2014StaticPDE} Fig. 3.  Although the biharmonic extension and linear extrapolation are of comparable accuracy in the vicinity of the interface, the biharmonic extension is free of the sharp corners that develop in the linear (and quadratic) extrapolation along the vertical centreline $x = 0$, which are artifacts of the cusp in the interface at this location.

\subsection{Example 2, three-dimensional peanut}
For our second example, we consider an extension of Example 1 to three dimensions over the domain $D = [-\pi, \pi]^3$.  Although the presentation to this point has concerned biharmonic extensions in two dimensions, all results generalise naturally to three.  The radial basis function associated with the extension is now the linear RBF $\Phi(r) = r$ (which is a fundamental solution to the biharmonic equation in three dimensions) and \eqref{eq:biharmonic} -- \eqref{eq:biharmonic_interpolation_constraints} continue to hold unchanged.  For this example we specify the field values
\begin{equation}
\label{eq:f_example1_3d}
f(x,y,z) = \cos(x)\sin(y)\sin(\pi/4-z)
\end{equation}
prescribed in the region
\[
\Omega = \left \{(x,y,z) \in D \left| \min\left(\sqrt{(x-0.8)^2 + y^2 + z^2} - 1, \sqrt{(x+0.8)^2 + y^2 + z^2} - 1\right) < 0 \right. \right \}\,.
\]
This gives rise to the peanut-shaped interface obtained by revolving the interface from Example 1 around the $x$ axis.  Again, the field is known on the inside and the goal is to smoothly extend the field to all of $D$.

This time we consider periodic conditions on the biharmonic extension.  In three dimensions, and especially with periodic conditions, the factorisation costs associated with the direct approach can be prohibitive owing to the amount of fill-in generated, so the iterative approach of Algorithm \ref{al:biharmonic} is beneficial (the necessary changes are to use symmetric array padding and discrete Fourier transforms in all three dimensions.)  Table \ref{tab:example_3d_errors} lists the runtimes, errors and estimated orders of convergence using this solution strategy.  The number of conjugate gradient iterations required to compute the field on an $n \times n \times n$ grid is seen to grow only slightly faster than linearly with $n$, requiring for example 112 and 276 iterations on $128^3$ and $256^3$ grids respectively.  By contrast, attempting to factorise the discrete biharmonic matrix quickly exhausted the available memory on the test machine.  Statistics for the constant, linear and quadratic extrapolations are also reported in Table \ref{tab:example_3d_errors}.  The biharmonic extension is observed to be several times faster to compute than the comparable, linear extrapolation.

\begin{table}
\centering
\caption{Runtimes, errors and estimated orders of convergence over a sequence of grids for Example 2 for the biharmonic extension (periodic conditions) and the constant, linear and quadratic extrapolations.  All errors are calculated with respect to the reference field \eqref{eq:f_example1_3d} in a neighbourhood of 4 grid cells of the interface.}
\begin{tabular}{l|*{4}{r}}
\hline
Grid & $32^3$ & $64^3$ & $128^3$ & $256^3$ \\
\hline
& \multicolumn{4}{c}{Biharmonic extension (Periodic)} \\
\hline
Iterations & 21 & 46 & 112 & 276 \\
Runtime & 0.1 & 1.5 & 37.6 & 733.9  \\
Error & 1.08e-01 & 3.42e-02 & 1.07e-02 & 2.84e-03  \\
Est. Order & -- & 1.66 & 1.68 & 1.91  \\
\hline
& \multicolumn{4}{c}{Constant extrapolation} \\
\hline
Runtime & 0.2 & 5.3 & 85.9 & 1269.9 \\
Error & 1.29e-01 & 6.87e-02 & 3.55e-02 & 1.77e-02  \\
Est. Order & -- & 0.91 & 0.95 & 1.00 \\
\hline
& \multicolumn{4}{c}{Linear extrapolation} \\
\hline
Runtime & 0.4 & 10.3 & 150.4 & 2258.8 \\
Error & 5.17e-02 & 1.14e-02 & 2.93e-03 & 7.42e-04  \\
Est. Order & -- & 2.18 & 1.96 & 1.98   \\
\hline
& \multicolumn{4}{c}{Quadratic extrapolation} \\
\hline
Runtime & 0.5 & 15.6 & 229.0 & 3420.3  \\
Error & 2.99e-02 & 3.31e-03 & 4.52e-04 & 7.66e-05  \\
Est. Order & -- & 3.18 & 2.87 & 2.56   \\
\hline
\end{tabular}
\label{tab:example_3d_errors}
\end{table}

\begin{figure}
\centering
\begin{tabular}{cc}
\subfloat[Reference field]{\includegraphics[trim= 2.5cm 0cm 2.5cm 0cm, clip, width=0.3\linewidth]{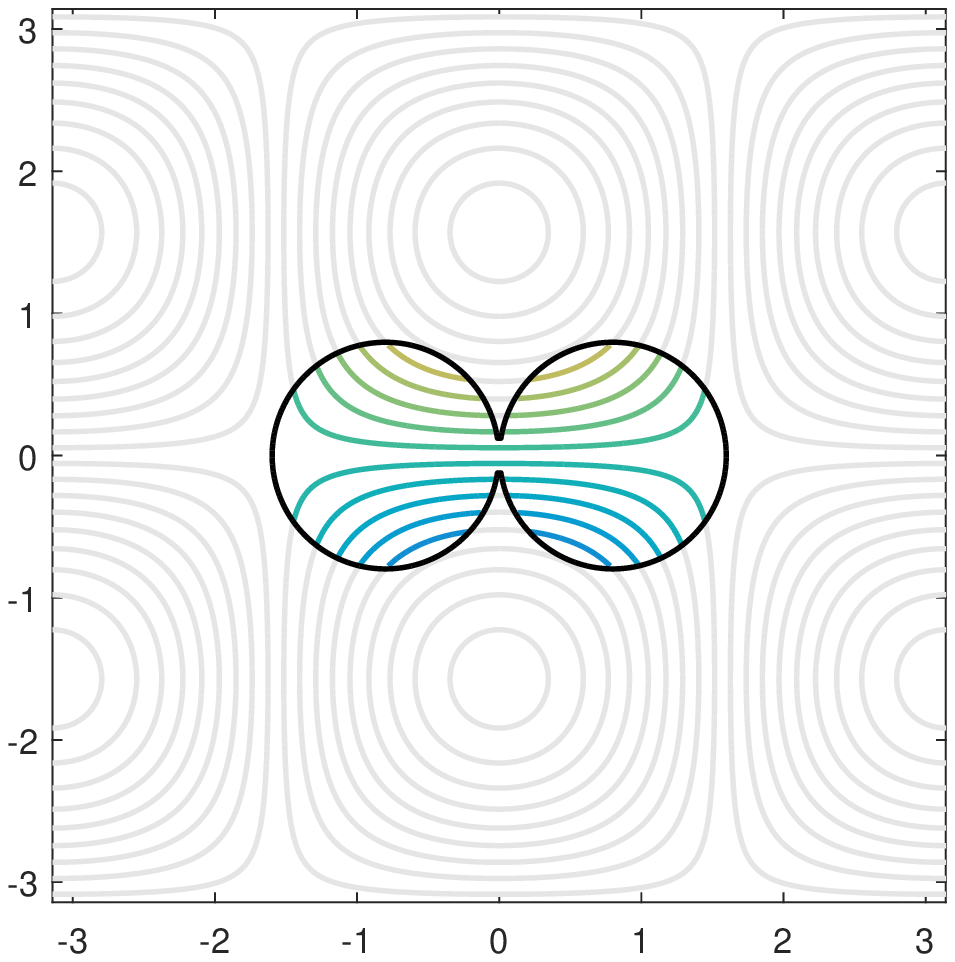}} &
\subfloat[Biharmonic extension (periodic)]{\includegraphics[trim= 2.5cm 0cm 2.5cm 0cm, clip, width=0.3\linewidth]{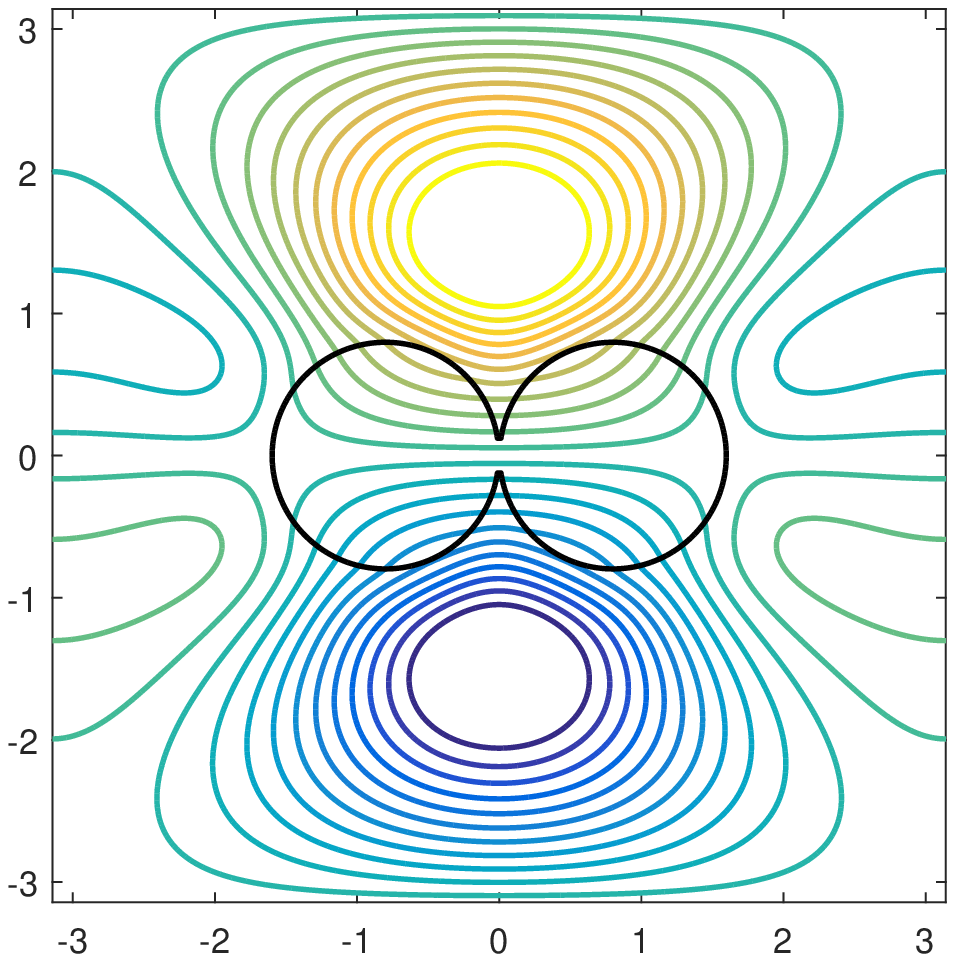}} \\
\subfloat[Constant extrapolation]{\includegraphics[trim= 2.5cm 0cm 2.5cm 0cm, clip, width=0.3\linewidth]{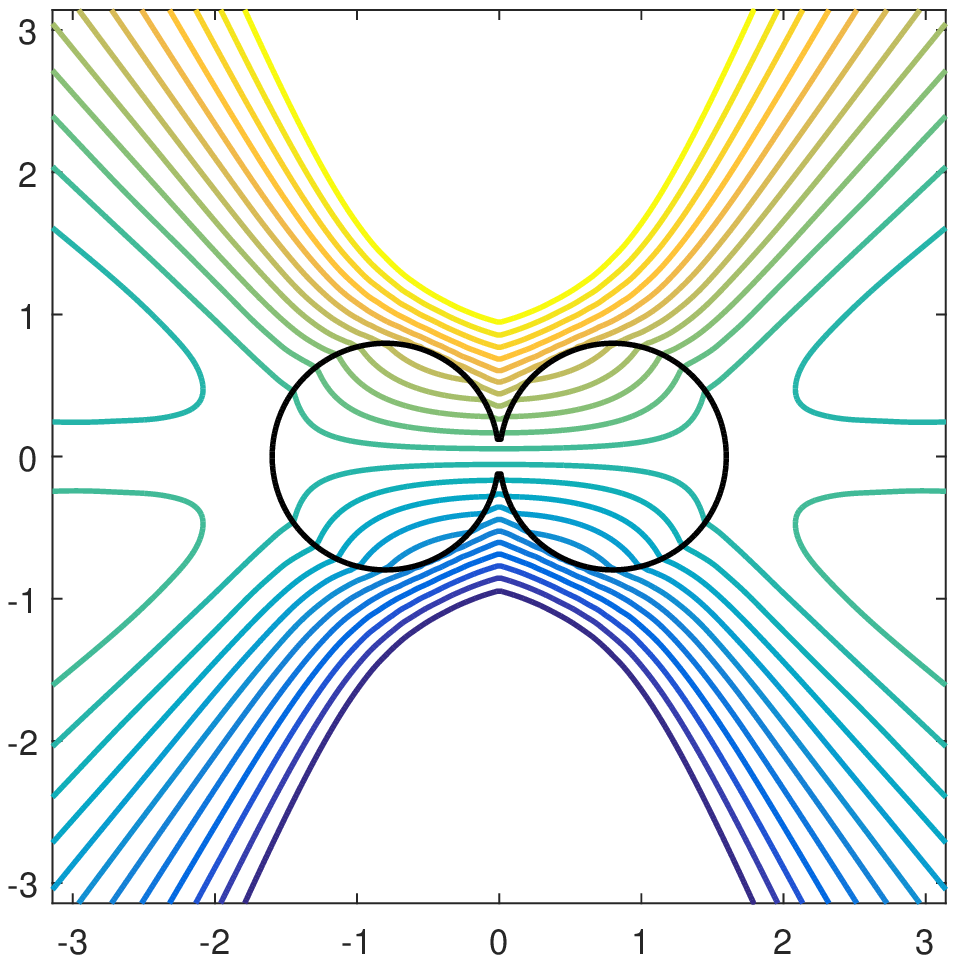}} &
\subfloat[Linear extrapolation]{\includegraphics[trim= 2.5cm 0cm 2.5cm 0cm, clip, width=0.3\linewidth]{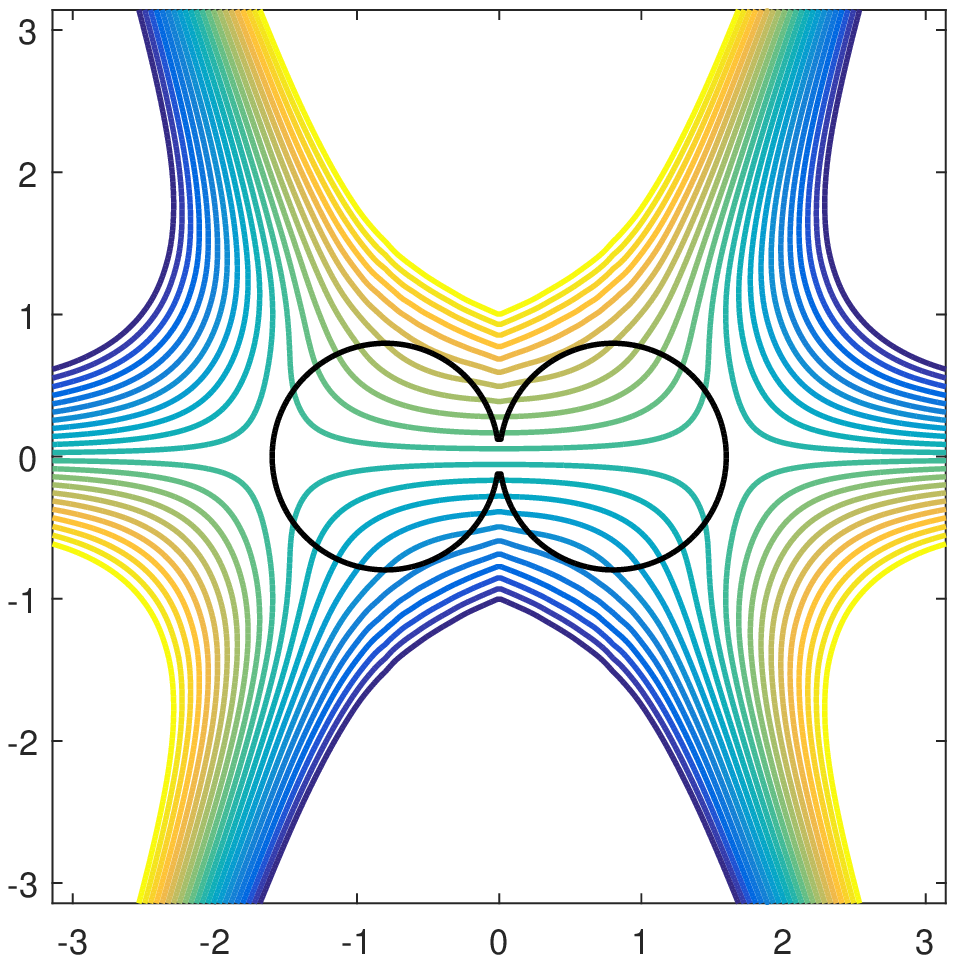}} \\
\subfloat[Quadratic extrapolation]{\includegraphics[trim= 2.5cm 0cm 2.5cm 0cm, clip, width=0.3\linewidth]{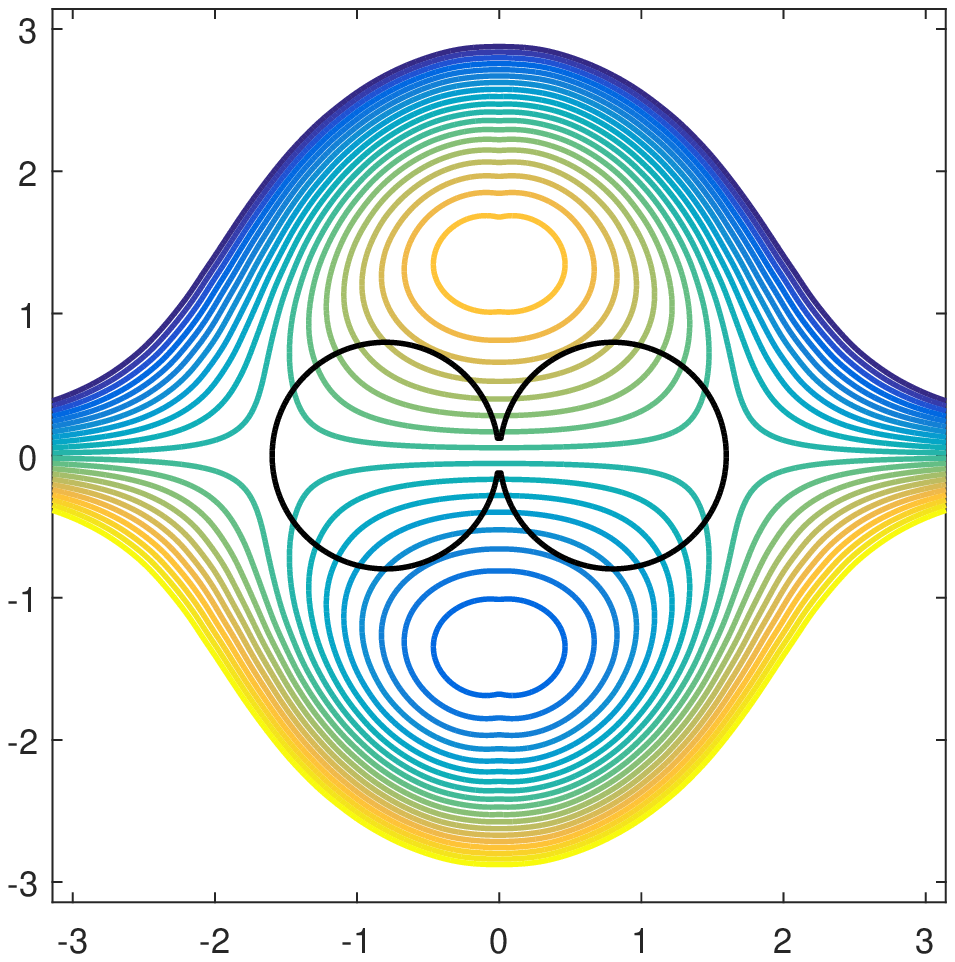}}
\end{tabular}
\caption{The reference field, periodic biharmonic extension, constant, linear and quadratic extrapolation for Example 2 at $z \approx 0.6$ on a $256 \times 256$ grid.  Contour levels are from $[-0.2, 0.2]$.}
\label{fig:example_3d}
\end{figure}

In Figure \ref{fig:example_3d} we present contours for the $z \approx 0.6$ slice for this problem.  All fields are accurate representations of the reference field in the vicinity of the interface.  Away from the interface, the biharmonic extension is free from the steep gradients that are present especially in the linear and quadratic extrapolations.

\subsection{Example 3, an annulus}
\begin{figure}
\centering
\begin{tabular}{cc}
\subfloat[Reference field]{\includegraphics[trim= 2.5cm 0cm 2.5cm 0cm, clip, width=0.45\linewidth]{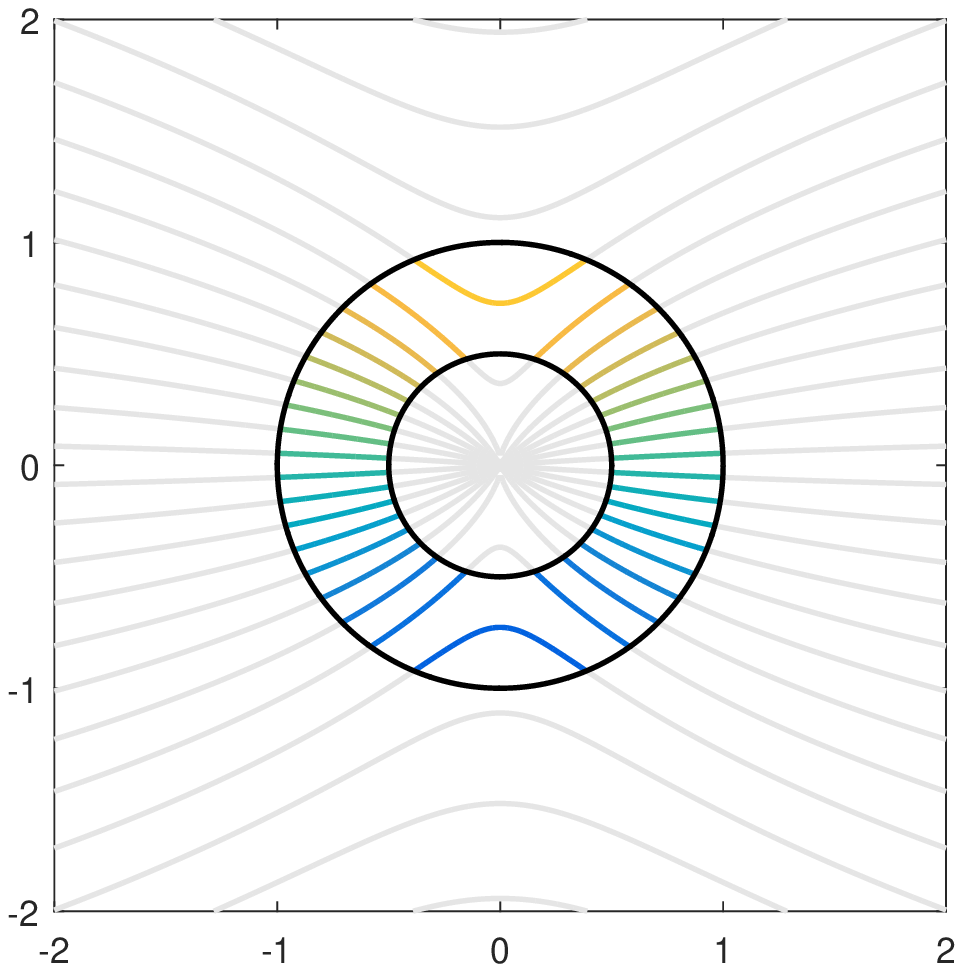}} &
\subfloat[Biharmonic extension subject to Dirichlet conditions]{\includegraphics[trim= 2.5cm 0cm 2.5cm 0cm, clip, width=0.45\linewidth]{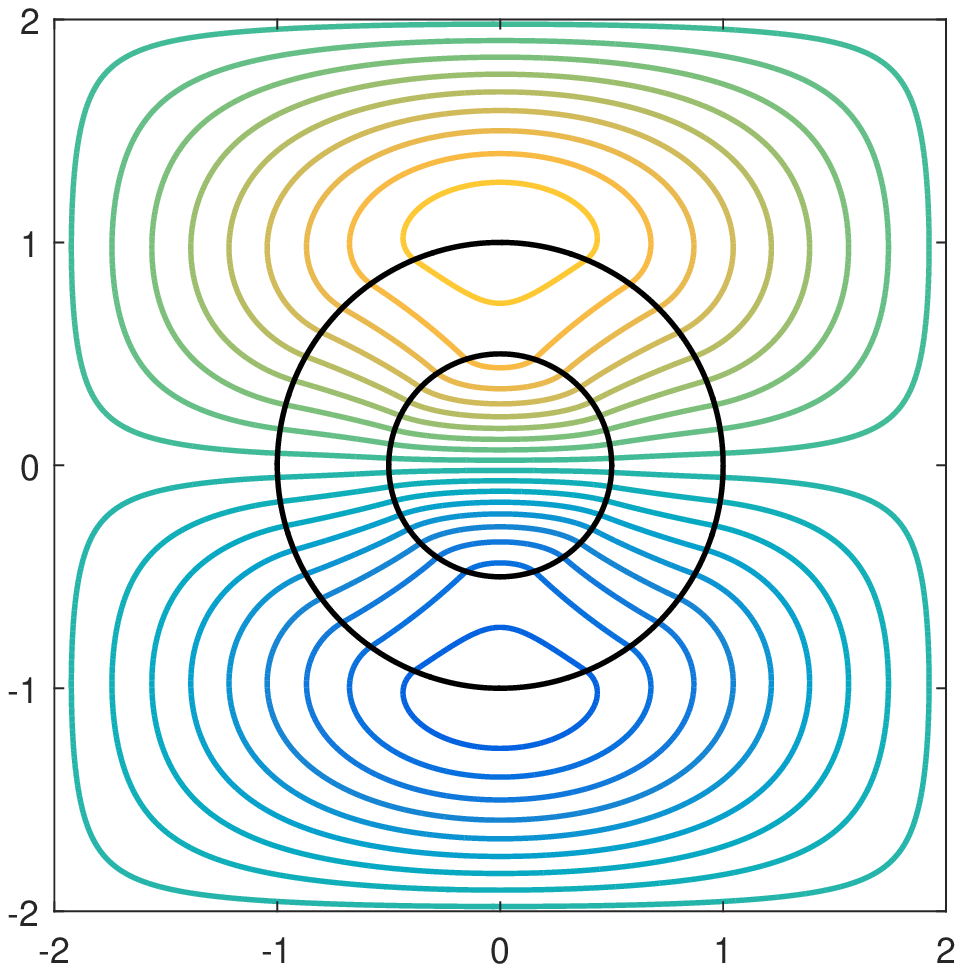}} \\
\subfloat[Biharmonic extension subject to Neumann conditions]{\includegraphics[trim= 2.5cm 0cm 2.5cm 0cm, clip, width=0.45\linewidth]{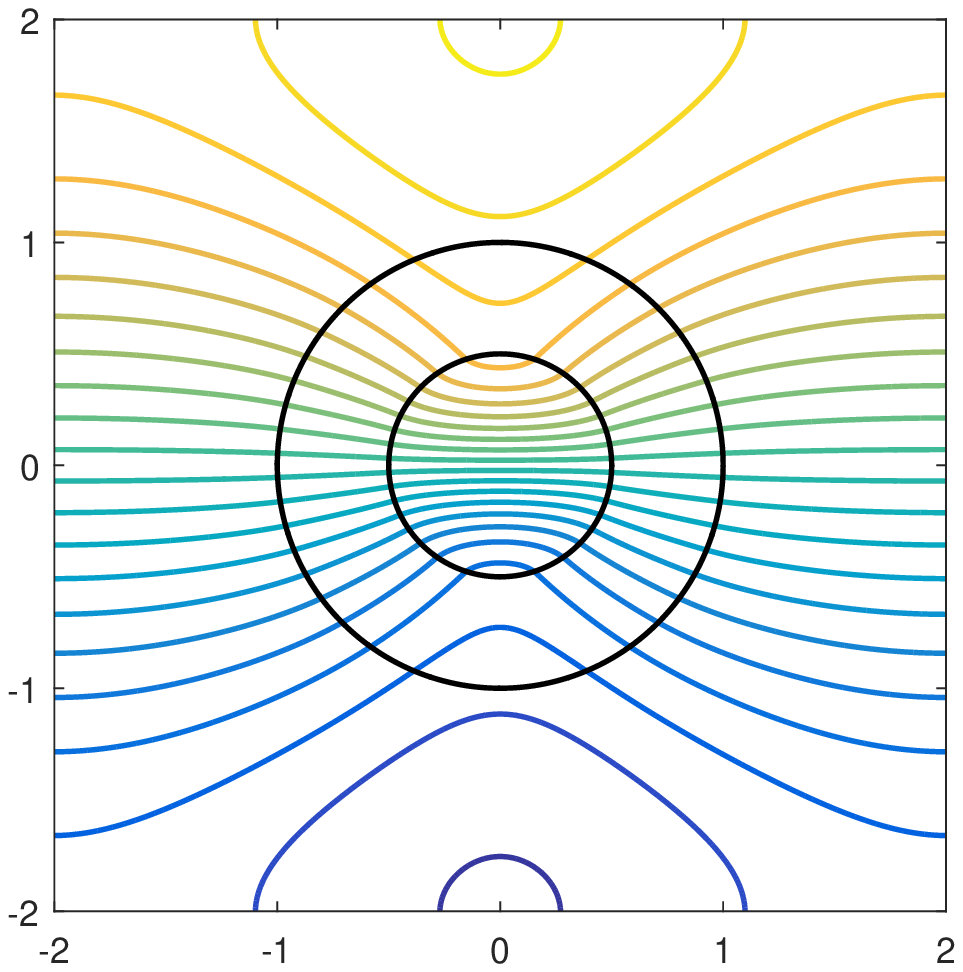}} &
\subfloat[Biharmonic extension subject to Neumann conditions and symmetry on a half grid (reference field also shown)]{\includegraphics[trim= 4.7cm 0cm 4.7cm 0cm, clip, width=0.23\linewidth, height=0.51\linewidth]{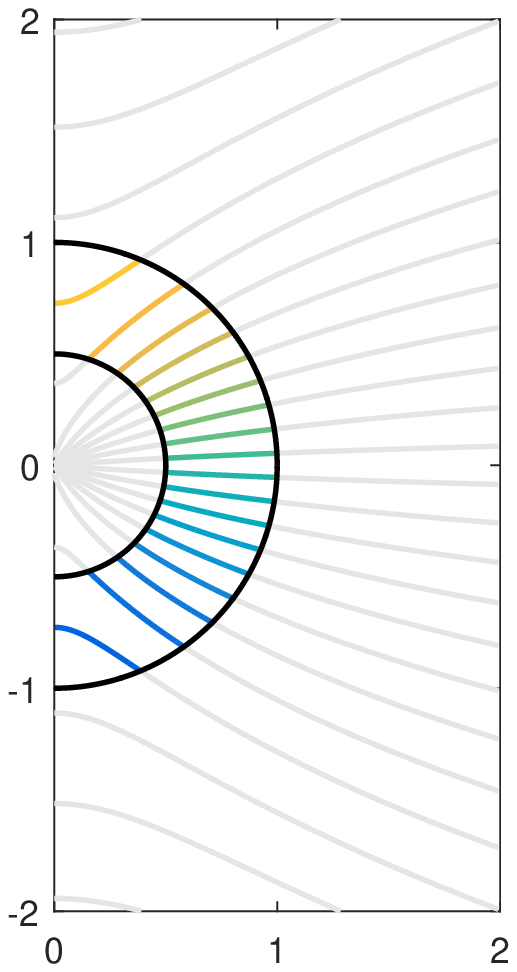}
\includegraphics[trim= 4.7cm 0cm 4.7cm 0cm, clip, width=0.23\linewidth, height=0.51\linewidth]{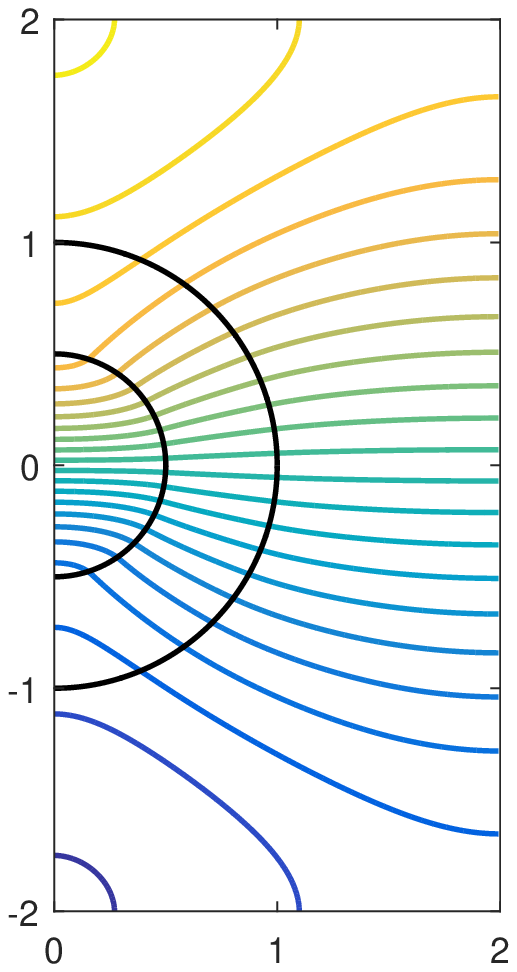}}
\end{tabular}
\caption{The reference field and biharmonic extension fields for Example 2 on a $256 \times 256$ grid subject to Dirichlet and Neumann boundary conditions, and on a half grid ($128 \times 256$) subject to Neumann conditions and symmetry. Contour levels are from $[-1.8, 1.8]$.}
\label{fig:example2}
\end{figure}

For our third example we take the domain $D = [-2, 2] \times [-2, 2]$ with field values
\begin{equation}
\label{eq:f_example2}
f(x,y) = y / \log(1 + \sqrt{x^2 + y^2})
\end{equation}
prescribed in the annulus
\[
\Omega = \left \{(x,y) \in D \left| 1/2 < \sqrt{x^2 + y^2} < 1 \right. \right \}\,.
\]
The two interfaces defining the annulus, and the known field within it, are illustrated in Figure \ref{fig:example2}(a).  Again, we also show in grey the contour lines for the reference field on all of $D$ for later comparison.  Multiply-connected domains such as this pose no additional difficulty to the biharmonic extension.  Indeed, the biharmonic extension is quite convenient in the sense that it extends the field to all of $D$ with a single solve (i.e. it extends both ``inside'' and ``outside'' at the same time).  We note that although the reference field is singular at the origin, the computed extensions are well defined and smooth at all points in $D$.

\begin{table}
\centering
\caption{Iteration counts, errors and estimated orders of convergence over a sequence of grids for Example 3.  Dirichlet and Neumann boundary conditions are tested, and errors are reported for inside and outside the annulus.  All errors are calculated with respect to the reference field \eqref{eq:f_example2} in a neighbourhood of 4 grid cells of the interface.}
\begin{tabular}{l|rrrr}
\hline
Grid& $128 \times 128$ & $256 \times 256$ & $512 \times 512$ & $1024 \times 1024$ \\
\hline
& \multicolumn{4}{c}{Dirichlet}\\
\hline
Iterations &         56 &        123 &        260 &        529 \\
Error outside &   6.15e-02 & 1.74e-02 & 4.66e-03 & 1.22e-03  \\
Est. Order outside &        -- &       1.83 & 1.90 & 1.94  \\
Error inside &   9.73e-02 &   2.80e-02 &   7.76e-03 &   2.07e-03 \\
Est. Order inside &        -- &       1.80 & 1.85 & 1.90 \\
\hline
& \multicolumn{4}{c}{Neumann}\\
\hline
Iterations &         57 &        124 &        256 &        522 \\
Error outside &   4.06e-03 & 1.07e-03 & 3.00e-04 & 7.91e-05  \\
Est. Order outside &        -- &       1.92 & 1.83 & 1.92 \\
Error inside &   9.73e-02 &   2.80e-02 &   7.76e-03 &   2.07e-03 \\
Est. Order inside &        -- &       1.80 & 1.85 & 1.90 \\
\hline
\end{tabular}
\label{tab:example2}
\end{table}

Figure \ref{fig:example2}(b)--(c) illustrates the biharmonic extension fields subject to Dirichlet and Neumann boundary conditions.  In the interior, the two extensions are identical, since the biharmonic stencils for interior nodes only reference other interior nodes or specified values within the annulus.  On the outside, both extensions appear similar near the interface, with the Neumann conditions producing an extension more in keeping with the reference field away from the interface.  Table \ref{tab:example2} confirms that the error within four grid cells on the outside is in fact smaller for the Neumann extension than for the Dirichlet, but both extensions converge at the expected rate to the reference field as the grid spacing is halved.  Once again, the number of conjugate gradient iterations grows slightly faster than linearly with $n$, up to a maximum of 522 iterations for $n = 1024$ with Neumann boundary conditions.

For this problem, an alternative approach is to exploit the symmetry of the interface and the field about the $y$ axis.  Applying Neumann conditions on the reduced problem over the rectangle $D = [0, 2] \times [-2, 2]$ produces the field shown in Figure \ref{fig:example2}(d): identical to that of Figure \ref{fig:example2}(c), but obtained using a grid with half as many nodes.  The number of conjugate gradient iterations required for the $512 \times 1024$ grid was also reduced compared to the full grid, down to 387 from 522.

\begin{figure}
\centering
\begin{tabular}{cc}
\subfloat[Biharmonic extension (Neumann)]{\includegraphics[trim= 2.5cm 0cm 2.5cm 0cm, clip, width=0.45\linewidth]{Example2_neumann_2.eps}} &
\subfloat[Constant extrapolation]{\includegraphics[trim= 2.5cm 0cm 2.5cm 0cm, clip, width=0.45\linewidth]{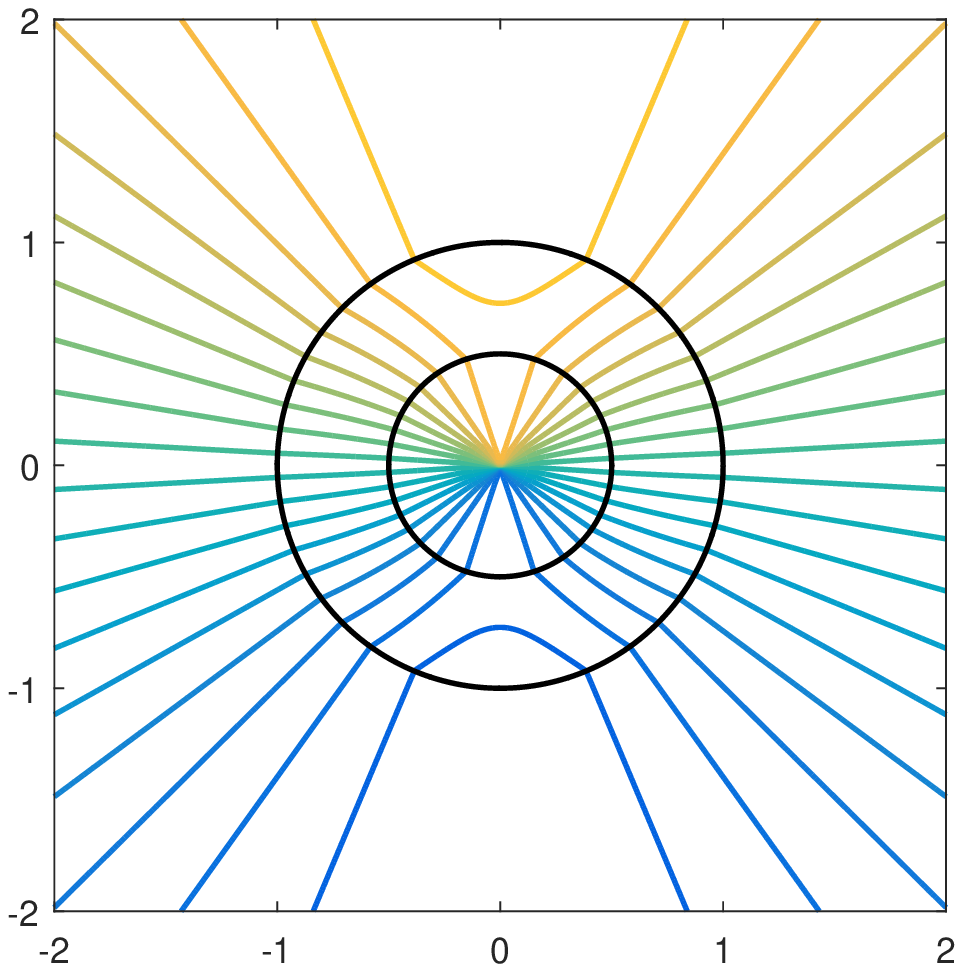}} \\
\subfloat[Linear extrapolation]{\includegraphics[trim= 2.5cm 0cm 2.5cm 0cm, clip, width=0.45\linewidth]{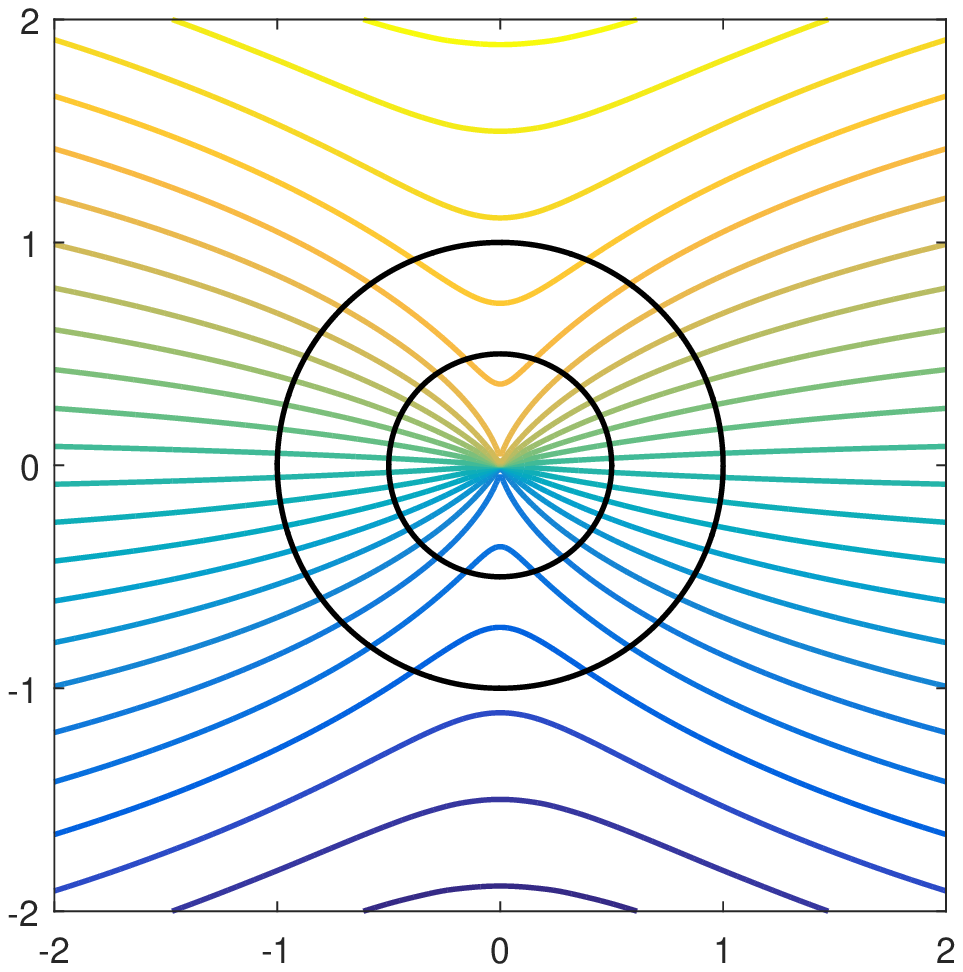}} &
\subfloat[Quadratic extrapolation]{\includegraphics[trim= 2.5cm 0cm 2.5cm 0cm, clip, width=0.45\linewidth]{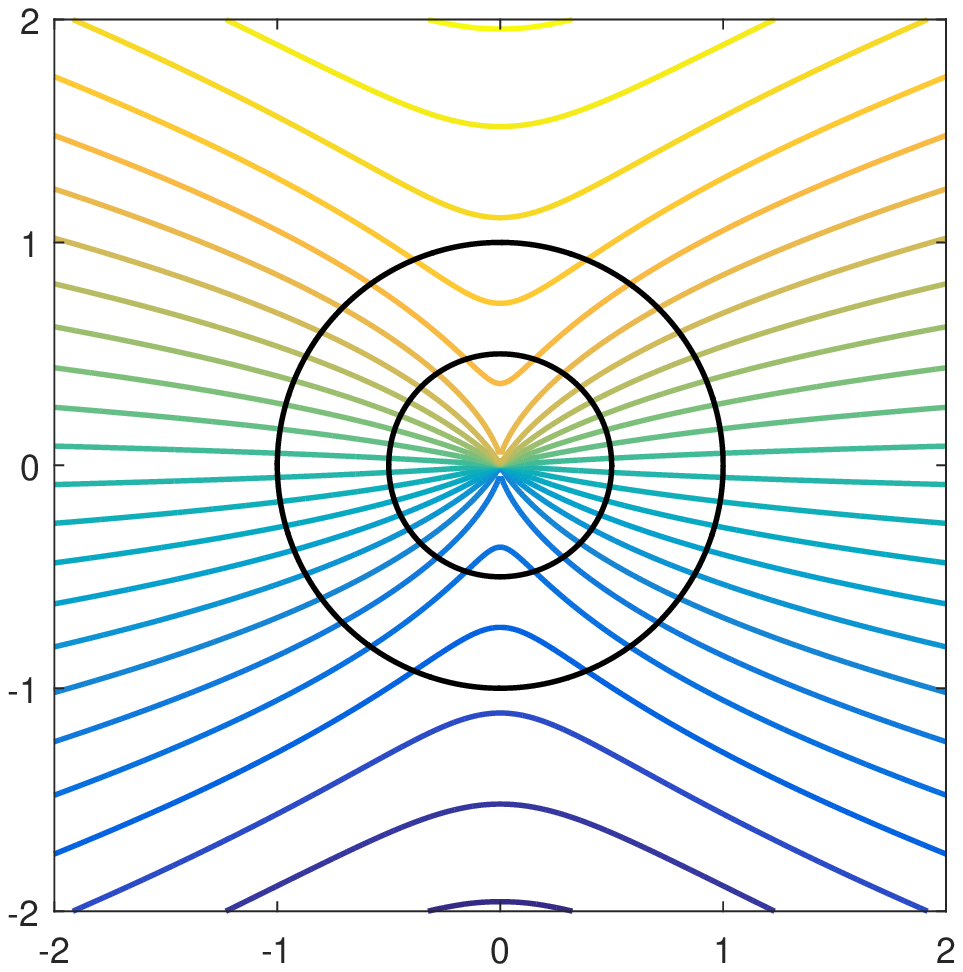}}\\
\end{tabular}
\caption{The biharmonic extension field for Example 3 on a $256 \times 256$ grid subject to Neumann boundary conditions, compared with constant extrapolation, linear extrapolation and quadratic extrapolation using the method of Aslam \cite{Aslam2004PartialDifferentialEquation}. Contour levels are from $[-1.8, 1.8]$.}
\label{fig:example2b}
\end{figure}

The biharmonic extension (with Neumann conditions) is compared with constant, linear and quadratic extrapolations in Figure \ref{fig:example2b}.  Here we note that these extrapolations do inherit the singularity at the origin from the field equation \eqref{eq:f_example2}, which the biharmonic extension does not.  Nonetheless in the vicinity of the interface all computed fields exhibit their expected orders of convergence, as confirmed in Table \ref{tab:example3_errors}.  In this table we also report runtimes for all methods.  For this two-dimensional problem, the direct method for the biharmonic extension is easily the fastest, at just 11.3 seconds, while the constant extrapolation and the iterative method for the biharmonic are of comparable runtime, at around 128 seconds.  The linear and quadratic extrapolations require, as expected, two and three times longer to compute than this, respectively.

\begin{table}
\centering
\caption{Runtimes, errors and ratios of errors over a sequence of grids for Example 2.  Dirichlet and Neumann boundary conditions are tested, and errors and ratios are reported for inside and outside the annulus.  All errors are calculated with respect to the reference field \eqref{eq:f_example2} in a neighbourhood of 4 grid cells of the interface.}
\begin{tabular}{l|*{4}{r}}
\hline
Grid & $128\times128$ & $256\times256$ & $512\times512$ & $1024\times1024$ \\
\hline
& \multicolumn{4}{c}{Biharmonic extension (Neumann)} \\
\hline
Runtime (direct) & 0.1 & 0.4 & 2.1 & 11.3  \\
Runtime (iterative) & 0.2 & 0.9 & 13.1 & 121.4 \\
Cond. Est. & 2.0e+07 & 3.1e+08 & 4.9e+09 & 7.8e+10 \\
Error outside &   4.06e-03 & 1.07e-03 & 3.00e-04 & 7.91e-05  \\
Est. Order outside &        -- &       1.92 & 1.83 & 1.92 \\
Error inside &   9.73e-02 &   2.80e-02 &   7.76e-03 &   2.07e-03 \\
Est. Order inside &        -- &       1.80 & 1.85 & 1.90 \\
\hline
& \multicolumn{4}{c}{Constant extrapolation} \\
\hline
Runtime & 0.2 & 1.0 & 12.9 & 128.3 \\
Error outside & 4.48e-02 & 2.33e-02 & 1.18e-02 & 6.05e-03 \\
Est. Order outside & -- & 0.94 & 0.99 & 0.96 \\
Error inside &   4.92e-02 & 2.45e-02 & 1.26e-02 & 6.75e-03  \\
Est. Order inside &        -- &       1.00 & 0.96 & 0.90 \\
\hline
& \multicolumn{4}{c}{Linear extrapolation} \\
\hline
Runtime & 0.3 & 2.0 & 25.7 & 259.7 \\
Error outside & 7.64e-04 & 1.77e-04 & 4.26e-05 & 1.06e-05 \\
Est. Order outside & -- & 2.11 & 2.06 & 2.00  \\
Error inside &   1.06e-03 & 2.27e-04 & 5.65e-05 & 1.62e-05  \\
Est. Order inside &        -- &       2.22 & 2.01 & 1.81 \\
\hline
& \multicolumn{4}{c}{Quadratic extrapolation} \\
\hline
Runtime & 0.5 & 3.0 & 39.4 & 389.5 \\
Error outside & 1.43e-04 & 1.83e-05 & 2.40e-06 & 3.05e-07 \\
Est. Order outside & -- & 2.96 & 2.93 & 2.97  \\
Error inside &   1.49e-03 & 1.90e-04 & 2.21e-05 & 2.77e-06  \\
Est. Order inside &        -- &       2.97 & 3.10 & 3.00 \\
\hline
\end{tabular}
\label{tab:example3_errors}
\end{table}

\subsection{Example 4, periodic domain / channel}
\begin{figure}
\centering
\begin{tabular}{ccc}
\subfloat[Reference field]{\includegraphics[trim= 2.5cm 0cm 2.5cm 0cm, clip, width=0.33\linewidth]{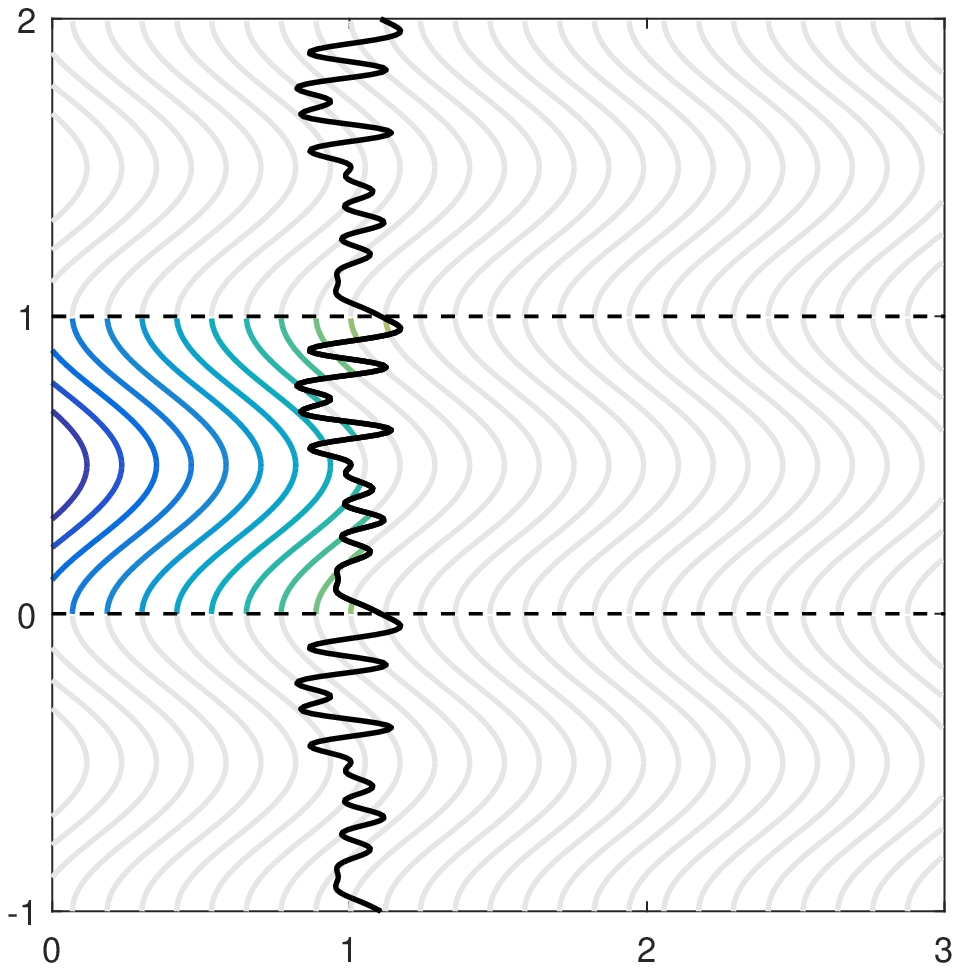}} &
\subfloat[Biharmonic extension subject to periodic and Neumann conditions]{\includegraphics[trim= 2.5cm 0cm 2.5cm 0cm, clip, width=0.33\linewidth]{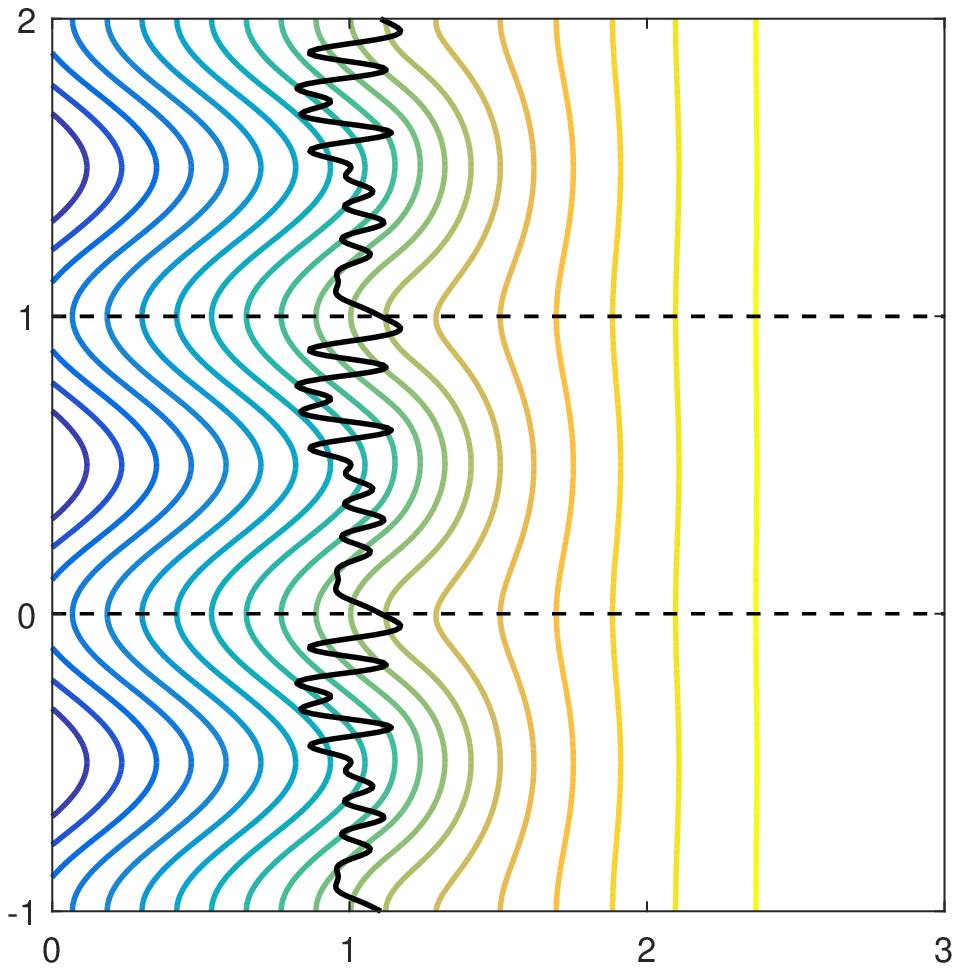}} &
\subfloat[Magnitude of error ($\log_{10}$) in a neighbourhood of the interface]{\includegraphics[trim = 4cm 0.3cm 4.5cm 0cm, clip, width=0.3\linewidth]{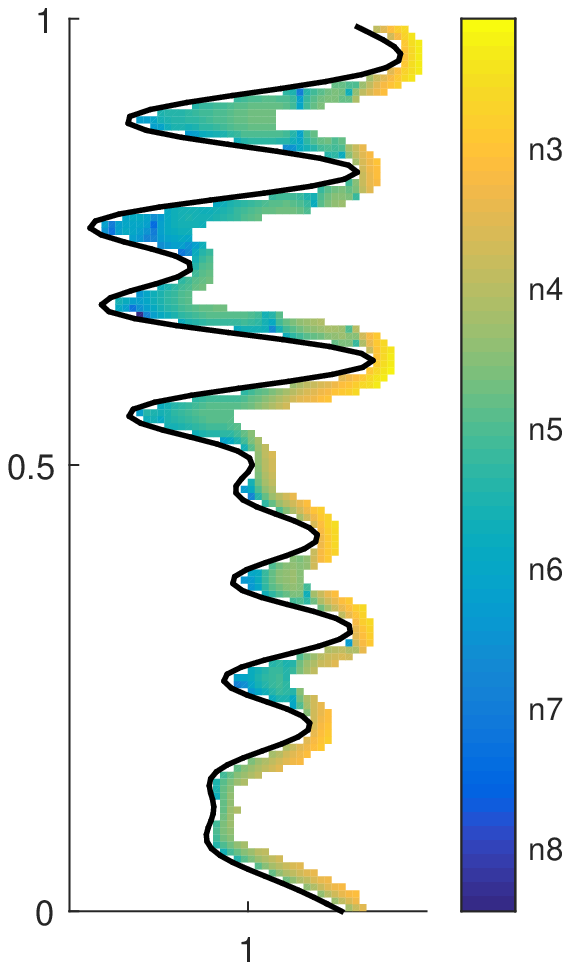}} \\
\end{tabular}
\caption{The extension field for Example 4 on a $384 \times 128$ grid subject to periodic and Neumann conditions.  Also shown is the magnitude of the error in a neighbourhood of four grid cells of the interface.}
\label{fig:example3}
\end{figure}

For our next example we take the domain $D = [0, 3] \times [0, 1]$ with periodic conditions on the top ($y = 1$) and bottom ($y = 0$).  The field values
\begin{equation}
\label{eq:f_example3}
f(x,y) = x+\cos(2\pi (y-0.25))/5
\end{equation}
are given on the left of the interface as shown in Figure \ref{fig:example3}(a).  The interface itself is described by a sum of sine waves $\sum_{k=1}^{10} a_k \sin(2\pi k y + \omega_k)$ with uniformly random amplitudes $a_k \in [-0.05, 0.05]$ and uniformly random phase shifts $\omega_k \in [0, 2\pi)$.  This problem is designed to simulate an application where the stability of a travelling front to small sinusoidal perturbations is investigated.  The periodic conditions are quite essential to the analysis, and hence the extension field must also respect this periodicity.  The biharmonic extension naturally does so with the imposition of periodic conditions on the top and bottom, and we chose Neumann conditions for the sides.

Figure \ref{fig:example3}(b) shows the biharmonic extension for this problem, replicated over two additional periods for illustrative purposes (the problem is solved numerically over just one period).  No special handling of nodes near the interface is required, as usual -- the biharmonic extension naturally takes care of incorporating the required field information whenever a stencil crosses the interface.  Close inspection of Figure \ref{fig:example3}(a) shows that some contour lines are cut off by the interface multiple times, alternately to the left and then to the right, but these and all other contours are smoothly extended in Figure \ref{fig:example3}(b).  However, we do observe in Table \ref{tab:example3} that the convergence rate of the biharmonic extension to the reference field is slower than in previous examples, estimated at only 1.5.

\begin{table}
\centering
\caption{Iteration counts, condition numbers, errors and estimated orders of convergence over a sequence of grids for Example 4 with periodic conditions on the top and bottom.  Two interfaces are tested: the irregular interface pictured in Figure \ref{fig:example3} consisting of 10 terms of scaled and shifted sine waves with increasing wavenumber, and a simpler interface consisting of just the first term.  All errors are calculated with respect to the reference field \eqref{eq:f_example3} in a neighbourhood of 4 grid cells of the interface.}
\begin{tabular}{l|rrrrrr}
\hline
Grid& $192 \times 64$ & $384 \times 128$ & $768 \times 256$ & $1536 \times 512$ \\
\hline
\hline
& \multicolumn{4}{c}{10 terms}\\
\hline
Iterations &         86 & 185 & 424 & 979 \\
Runtime (iterative) & 0.3 & 1.1 & 19.6 & 200.3 \\
Runtime (direct) & 0.0 & 0.2 & 1.4 & 7.0  \\
Cond. Est. & 5.8e+08 & 9.1e+09 & 1.4e+11 & 2.3e+12 \\
Error &   2.45e-02 & 9.27e-03 & 2.96e-03 & 1.03e-03 \\
Est. Order &        -- &       1.40 & 1.65 & 1.52 \\
\hline
& \multicolumn{4}{c}{1 term}\\
\hline
Iterations &         38 & 70 & 131 & 269 \\
Runtime (iterative) & 0.5 & 0.5 & 5.9 & 52.3 \\
Error &   1.54e-02 & 4.66e-03 & 1.30e-03 & 3.52e-04 \\
Error Order &        -- &       1.72 & 1.84 & 1.88 \\
\hline
\end{tabular}
\label{tab:example3}
\end{table}

The lower convergence rate is due to the irregularity of the interface, resulting in biharmonic stencils near the peaks of the interface incorporating fewer nodes on the inside that carry known field values.  In Figure \ref{fig:example3}(c), we see that the largest errors between the reference and extension fields occur in these regions near the peaks.  Conversely, the errors are much lower in regions near the troughs of the interface, where the stencils reference many field values on the inside.  If, instead, we use only a single period of a sine wave $a_1 \sin(2\pi y + \omega_1)$ as the shape of the interface (not pictured), we obtain the second set of results in Table \ref{tab:example3}.  Here, convergence is restored to near second order, confirming that the irregularity of interface is the cause of the reduced convergence rate.

The preconditioned conjugate gradient iteration continues to be effective, although for the irregular interface on the finest grid with $N = 1536 \times 512$ nodes, the number of iterations required increased to more than $\sqrt{N}$, compared to less than $\sqrt{N}/2$ for the single period sine wave.  This suggests that as grids are further refined, more advanced iterative solvers may be worth investigating.  The direct solver is, naturally, unaffected by such concerns, and generates solutions more than an order of magnitude faster than the iterative solver on this two-dimensional problem.

\paragraph{Channel geometry}
Alternatively, the problem may be considered in a channel, with Neumann conditions imposed on the top and bottom.  For comparative purposes, we use the same interface, but now the prescribed field
\begin{equation}
\label{eq:f_example3b}
f(x,y) = x+\cos(\pi y)/5
\end{equation}
satisfies Neumann, rather than periodic, conditions at the top and bottom.

\begin{figure}
\centering
\begin{tabular}{cc}
\subfloat[Reference field]{\includegraphics[trim= 1.5cm 3cm 1.0cm 2.5cm, clip, width=0.5\linewidth]{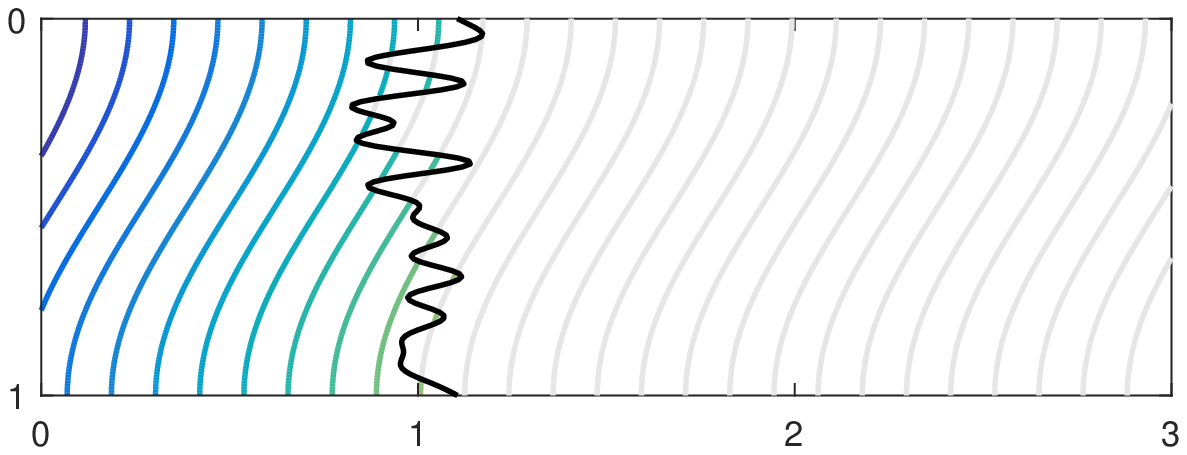}}
& \renewcommand{\thesubfigure}{c} \multirow{-7}[-4]{*}{
\subfloat[Magnitude of error ($\log_{10}$) in a neighbourhood of the interface]{\includegraphics[trim = 4cm 0.8cm 4.5cm 0cm, clip, width=0.3\linewidth]{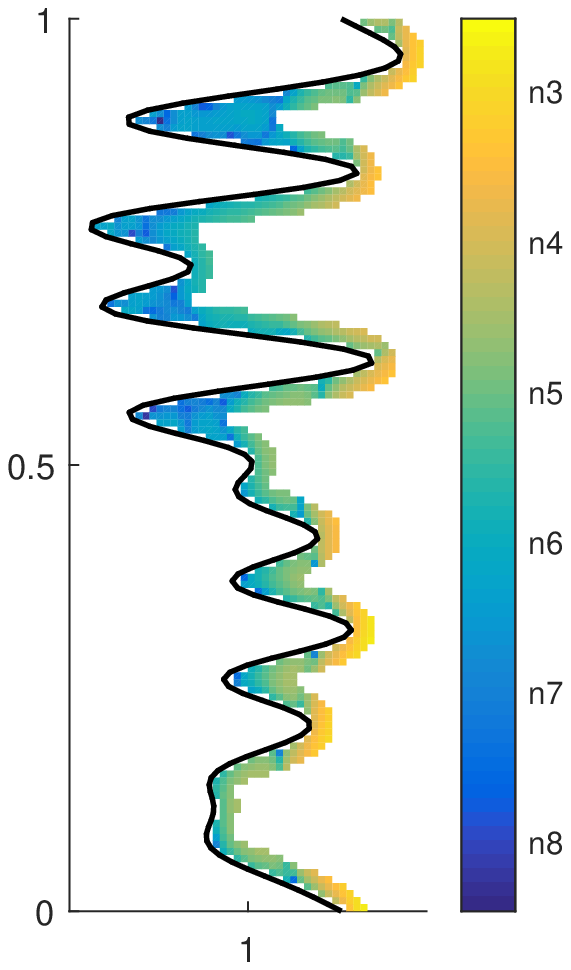}}} \\
\renewcommand{\thesubfigure}{b}
\subfloat[Biharmonic extension subject to Neumann conditions]{\includegraphics[trim= 1.5cm 3cm 1.0cm 2.5cm, clip, width=0.5\linewidth]{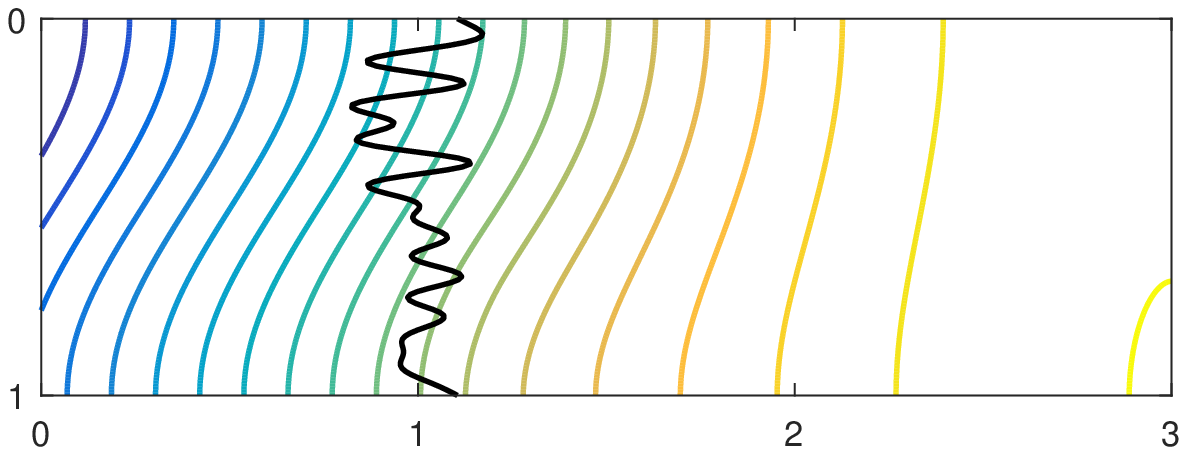}}\\
\end{tabular}
\caption{The extension field for Example 4 on a $384 \times 128$ grid subject to Neumann conditions.  Also shown is the magnitude of the error in a neighbourhood of four grid cells of the interface.}
\label{fig:example3b}
\end{figure}

We present our results in Figure \ref{fig:example3b} and Table \ref{tab:example3b}.  There are no significant differences observed in the behaviours of the runtimes or the errors in Table \ref{tab:example3b} as compared to Table \ref{tab:example3}.  From Figure \ref{fig:example3b}, we observe that the field is smoothly extended across the interface, and is in keeping with the Neumann conditions prescribed on the boundaries.  We note again that the error in the extended field is greatest near the peaks of the interface, due to the small number of neighbouring nodes that carry known field values.  Indeed, for the coarsest grid, $192 \times 64$, the biharmonic stencils at the tips of the steepest peaks reference just two values on the inside: $f_{i-1,j}$ and $f_{i-2,j}$, with all other neighbouring values being themselves constructed as part of the extension.  Although this scenario leads to a local reduction in the accuracy, it causes no particular problems for the biharmonic extension in its implementation.  In fact, it is worthwhile emphasising that even a single isolated node surrounded by the interface, or a collection of such isolated nodes, lead to well-posed biharmonic extensions.  In the latter case the method would simply be building a clamped thin plate spline that interpolates the data at the isolated points, which seems as reasonable an extension as one could hope for under those circumstances.

\begin{table}
\centering
\caption{Iteration counts, condition numbers, errors and estimated orders of convergence over a sequence of grids for Example 4 with Neumann conditions on the top and bottom.  Two interfaces are tested: the irregular interface pictured in Figure \ref{fig:example3b} consisting of 10 terms of scaled and shifted sine waves with increasing wavenumber, and a simpler interface consisting of just the first term.  All errors are calculated with respect to the reference field \eqref{eq:f_example3b} in a neighbourhood of 4 grid cells of the interface.}
\begin{tabular}{l|rrrrrr}
\hline
Grid& $192 \times 64$ & $384 \times 128$ & $768 \times 256$ & $1536 \times 512$ \\
\hline
\hline
& \multicolumn{4}{c}{10 terms}\\
\hline
Iterations &         97 & 212 & 496 & 1176 \\
Runtime (iterative) & 0.3 & 1.3 & 17.9 & 206.6 \\
Runtime (direct) & 0.1 & 0.2 & 1.1 & 6.3  \\
Cond. Est. & 5.6e+08 & 9.1e+09 & 1.5e+11 & 2.3e+12 \\
Error &   8.82e-03 & 3.11e-03 & 1.11e-03 & 3.76e-04 \\
Est. Order &        -- &       1.50 & 1.48 & 1.57  \\
\hline
& \multicolumn{4}{c}{1 term}\\
\hline
Iterations &         43 & 88 & 180 & 384 \\
Runtime (iterative) & 0.1 & 0.7 & 7.2 & 68.0 \\
Error &   5.88e-03 & 1.63e-03 & 4.54e-04 & 1.32e-04 \\
Error Order &        -- &       1.85 & 1.84 & 1.79 \\
\hline
\end{tabular}
\label{tab:example3b}
\end{table}

\subsection{Example 5, Stefan problem}
For our final example, we implement the biharmonic extension in the context of a full level set method to solve a Stefan problem.  The problem models solidification of a seed in an undercooled liquid, and is governed by the following dimensionless equations.  Temperature in both solid ($s$) and liquid ($\ell$) phases obeys the heat equation,
\begin{equation} \label{eq:sp_temp}
\frac{\partial T}{\partial t} = \nabla^2T \quad \text{in} \quad D.
\end{equation}
The normal velocity of the interface that separates the two phases is given by the jump condition
\begin{equation} \label{eq:sp_vel}
V_n = - \left( \frac{\partial T_\ell}{\partial n} - \frac{\partial T_s}{\partial n} \right) = - \left[ \frac{\partial T}{\partial n} \right] \quad \text{on} \quad \partial \Omega.
\end{equation}
The temperature on the interface is given by the Gibbs-Thomson relation,
\begin{equation} \label{eq:sp_int}
T = - \sigma \kappa \quad \text{on} \quad \partial \Omega,
\end{equation}
which models the effects of surface tension with coefficient $\sigma$, where $\kappa$ is the interface curvature.  The model is subject to the following initial condition for temperature,
\begin{equation} \label{eq:sp_ic}
T \left( \mathbf{x}, 0 \right) = \begin{dcases*} 0, & $\mathbf{x} \in \Omega$ \\ - \frac{1}{\beta}, & $\mathbf{x} \in D \backslash \Omega$ \end{dcases*},
\end{equation}
where $\beta$ is the Stefan number.
The temperature in the far-field is given by $T \rightarrow - 1 / \beta$ as $\|\mathbf{x}\| \rightarrow \infty$, which is imposed as Neumann boundary conditions on the computational domain $D$.

This problem was considered by Gibou et al. \cite{Gibou2013} who presented solutions based on a ghost-fluid level-set method previously proposed in \cite{gibou2003level}.  Briefly, with the ghost-fluid method, the heat equation \eqref{eq:sp_temp} is discretised using finite differences and, wherever computational stencils cross the interface, ghost values are constructed by extrapolation which implicitly satisfy the jump condition.  Temperate gradients are calculated in each phase, and the normal derivative is then extrapolated in the normal direction, allowing a velocity field to be calculated at all grid points using \eqref{eq:sp_vel}.  The level set method is then used to advance the moving front according to \eqref{eq:level_set_equation}, with the level set function reinitialised to approximate a signed distance function after each step.  Further details can be found in \cite{gibou2003level, Gibou2013, fedkiw1999non}.

Our numerical scheme does not make use of ghost cells, and instead uses asymmetric stencils in the vicinity of the interface in the manner described by Chen et al.~\cite{chen1997simple}. 
Hence we solve for only one temperature field at each node, that being the temperature of the phase present there.  We compute gradients from these temperatures, again using asymmetric stencils in the vicinity of the interface, in both the $x$-$y$ and $\eta$-$\zeta$ (i.e. $45^\circ$) directions, as suggested by Chen et al. \cite{chen1997simple}.  These gradient components, defined only at grid nodes, are then extended from one phase to the other using biharmonic extensions with Neumann boundary conditions.  That is, the solid phase gradients are extended outwards using one biharmonic solve, and the liquid phase gradients are extended inwards using another biharmonic solve, allowing a velocity field to be calculated at all grid points using \eqref{eq:sp_vel}.  From there the level set method proceeds as outlined above.

\psfrag{AA}[cc]{-1.5}
\psfrag{BB}[cc]{-0.5}
\psfrag{CC}[cc]{0.5}
\psfrag{DD}[cc]{1.5}

\begin{figure}
\centering
\begin{tabular}{ccc}
\subfloat[$\sigma = 0$]{\includegraphics[trim= 1.8cm 0.5cm 1.8cm 0.8cm, clip,width=0.32\linewidth]{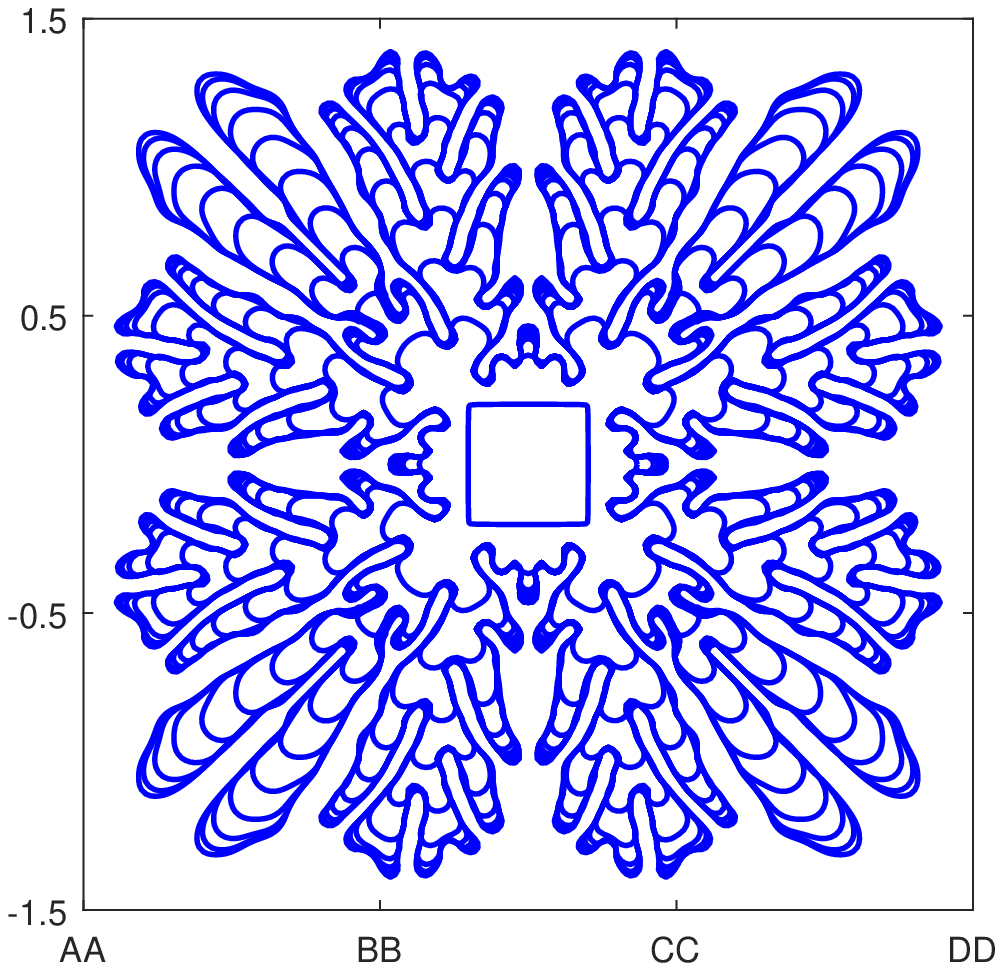}}
&
\subfloat[$\sigma = 0.0005$]{\includegraphics[trim= 1.8cm 0.5cm 1.8cm 0.8cm, clip,width=0.32\linewidth]{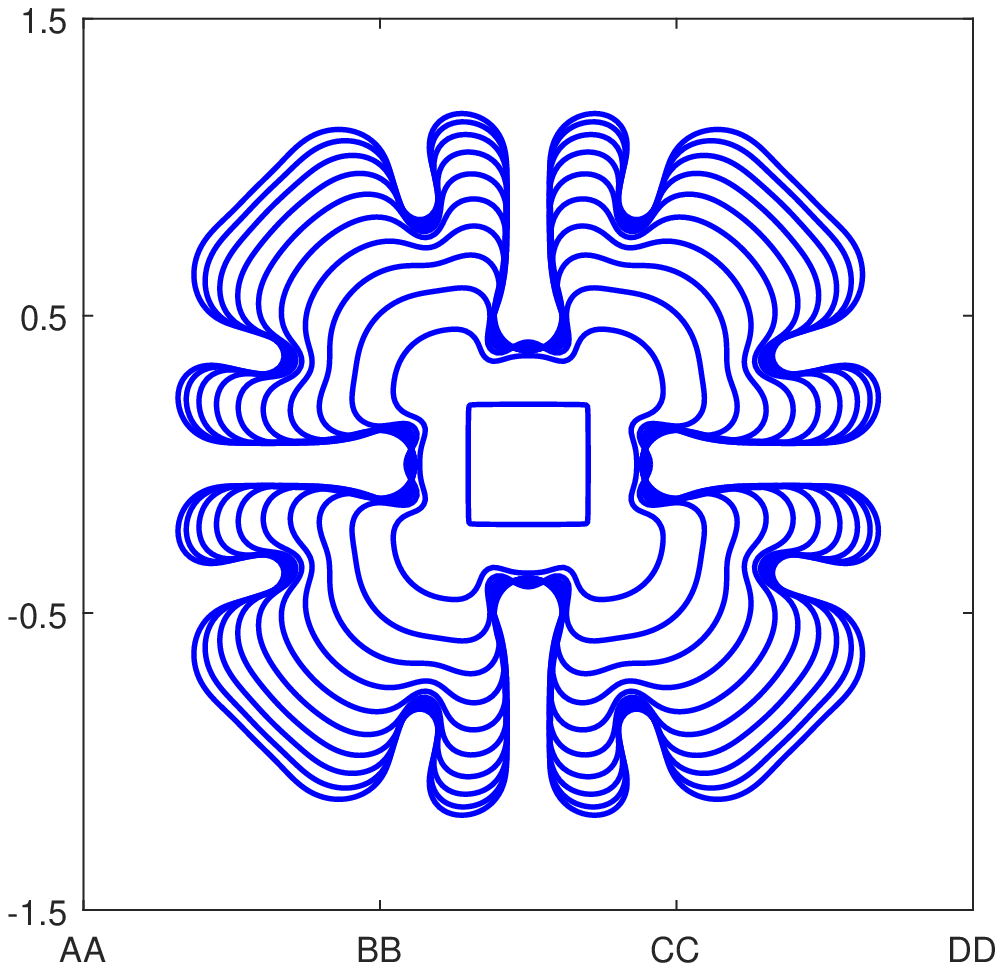}}
&
\subfloat[$\sigma = 0.001$]{\includegraphics[trim= 1.8cm 0.5cm 1.8cm 0.8cm, clip,width=0.32\linewidth]{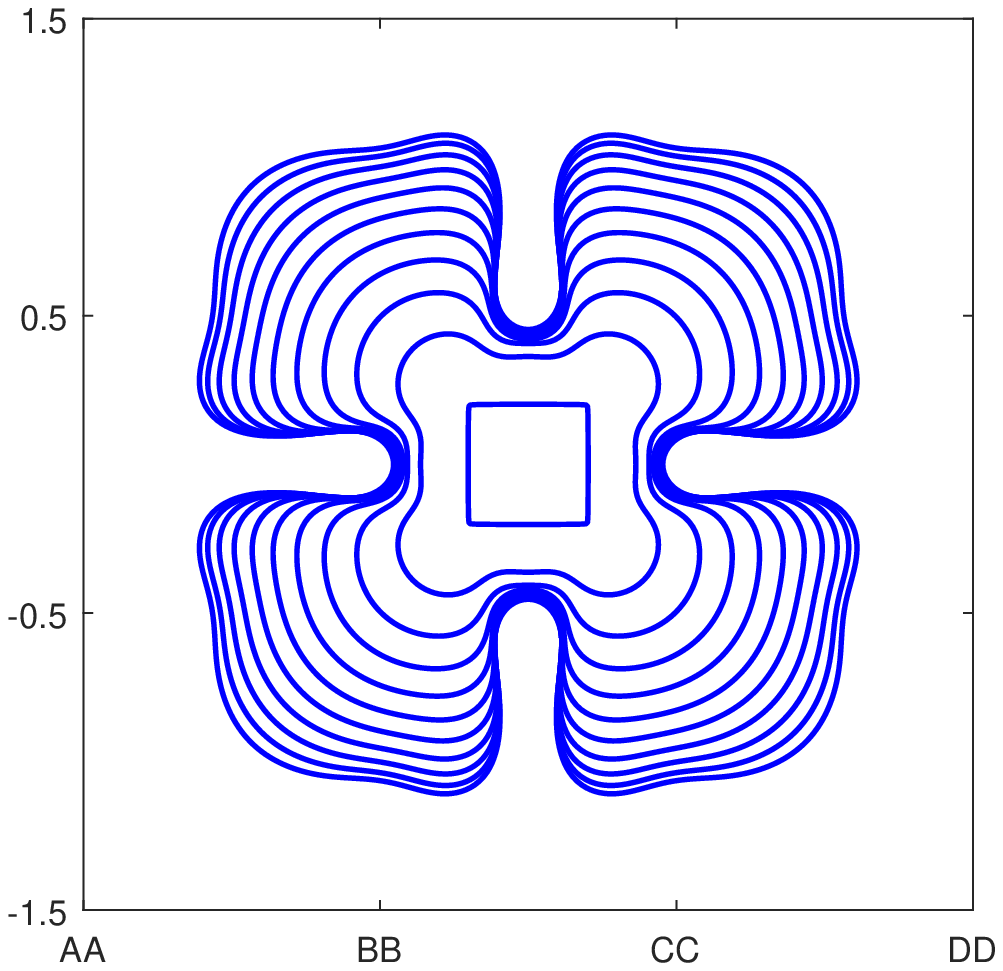}}\\
\end{tabular}
\caption{Solutions of the Stefan problem \eqref{eq:sp_temp} - \eqref{eq:sp_ic} to $t = 0.4$ on a $200 \times 200$ grid with $\beta = 2$ and varying $\sigma$.}
\label{fig:stefan}
\end{figure}

Figure \ref{fig:stefan} reproduces Fig. 24 of \cite{Gibou2013} using our approach.  The Stefan number is held fixed at $\beta = 2$, time runs to $t = 0.4$, and the surface tension coefficient $\sigma$ is varied.  Our solutions are shown for a $200 \times 200$ grid.  We used first order Euler time stepping schemes (implicit for the heat equation, explicit for the level set equation) with small fixed stepsize $\Delta t = 0.0005$.  We note that the Stefan problem with $\sigma = 0$ is ill-posed, and the numerical scheme is in effect providing a form of regularisation, which allows Figure \ref{fig:stefan}(a) to be obtained.  Hence, with different numerical schemes it is not expected to obtain agreement for this figure, and indeed we find our Figure \ref{fig:stefan}(a) is not identical to that published in \cite{Gibou2013}, although both solutions are highly unstable.  With $\sigma > 0$ the problem is well-posed (although still unstable), and our numerical scheme incorporating the biharmonic extension produces results (Figure \ref{fig:stefan}(b)-(c)) consistent with those published in \cite{Gibou2013}.

These results serve to illustrate the utility of the biharmonic extension as part of a level set solution strategy.  In our implementation, we utilised sparse direct factorisation in solving both the heat equation and the biharmonic equation.  We found that factorising the discrete biharmonic matrix took two to three times as long as factorising the (scaled, shifted) discrete Laplacian matrix of the heat equation.  This led to approximately 40\% of the total runtime being taken up with biharmonic solves, 17\% with heat equation solves, and the remaining runtime split between gradient calculations, level set reinitialisations, and other level set machinery.

\section{Conclusion}
\label{sec:conclusion}
In this paper we have presented a method for computing field extensions in level set methods by solving a biharmonic equation.  This approach differs from other typical approaches in its fully implicit treatment of the interface.  No explicit properties of the interface such as its location or the field value on the interface are required.  This leads to a particularly simple implementation using either a sparse factorisation, or a matrix-free conjugate gradient solver with a fast Poisson preconditioner.

We demonstrated the biharmonic extension on a number of test problems that served to illustrate its effectiveness at producing smooth and accurate extensions near interfaces, whether extending inside or outside, on single or multiply-connected domains, and conforming to a variety of prescribed boundary conditions including symmetry and periodicity.  As part of a level set solution strategy for moving interface problems, the method works well, as was demonstrated for an unstable Stefan problem.  Further, we have shown there is no essential difference in applying this method in three dimensions, given that the biharmonic equation yields the natural choice of radial basis function $\Phi(r) = r$ in three dimensions, as it does with $\Phi(r) = r^2 \log(r)$ in two.

Overall, we believe this approach will prove to be a very convenient method for incorporating velocity or other field extensions into level set methods.  Developing further improvements to the efficiency of the iterative solution strategy, particularly important in three dimensions, seems an excellent avenue for future research.

\section*{Acknowledgements}
TJM acknowledges the helpful discussions with Stephen Roberts from Australian National University.  SWM acknowledges the support of the Australian Research Council via the Discovery Project DP140100933.  All authors thank the anonymous referees for their suggestions that have helped improve this paper.

\end{document}